\documentstyle[11pt,leqno,amscd,amssymb,pstricks,latexsym,amsbsy,xypic,mathrsfs,verbatim]{amsart}

\input xy
\xyoption{all}

\newtheorem{Thm}{Theorem}[section]
\newtheorem{Lem}[Thm]{Lemma}

\newtheorem{Prop}[Thm]{Proposition}

 \theoremstyle{definition}

\newtheorem{Rem}[Thm]{Remark}

\numberwithin{equation}{section}

\newcommand\Coker{{\rm Coker}}

\renewcommand\mod{{\rm mod\,}}

\newcommand\Hom{{\rm Hom}}
\newcommand\End{{\rm End}}
\newcommand\Ext{{\rm Ext}}

\newcommand\Ker{{\rm Ker}}

\newcommand\coh{{\rm coh\,}}
\newcommand\vect{{\rm vect}}

\newcommand\bbZ{{\mathbb Z}}

\newcommand\bbY{{\mathbb Y}}
\newcommand\bbX{{\mathbb X}}

\newcommand\cS{{\mathcal S}}

\renewcommand\cD{{\mathcal D}}
\renewcommand\cH{{\mathcal H}}

\newcommand\cC{{\mathcal C}}
\newcommand\cU{{\mathcal U}}
\newcommand\cc{{\mathcal C}}

\newcommand{\blambda}{{\boldsymbol{\lambda}}}

\newcommand\vez{\varepsilon}

\newcommand{\id}{{\rm id}}
\newcommand{\add}{{\rm add}}

\newcommand\im{{\rm Im}}

\begin{document}

\title[Recollements and Ladders for weighted projective lines]
{Recollements and Ladders \\ for
weighted projective lines}

\author{Shiquan Ruan}
\address{School of Mathematical Sciences, Xiamen University,
Xiamen 361005, China}
\email{sqruan@@xmu.edu.cn}
\address{Faculty of Mathematics, Bielefeld University,
Bielefeld D-33501, Germany.} \email{sruan@@math.uni-bielefeld.de}

\begin{abstract} In this paper, we construct recollements and ladders for exceptional curves by using reduction/insertion functors due to $p$-cycle construction. As applications to weighted projective lines, we classify recollements for the category of coherent sheaves over a weighted projective line, and give an explicit description of ladders in two different levels: the bounded derived category of coherent sheaves and the stable category of vector bundles.
 \end{abstract}

%\date{\today}

\subjclass[2010]{14H45, 16E35, 16G20, 17B37}

\keywords{Weighted projective line, $p$-cycle, recollement, ladder, exceptional curve}

\maketitle

\section{Introduction}

This work is devoted to weighted projective lines and their categories of coherent sheaves. These arise frequently in representation theory \cite{[GL1],[L],[L1]}. Each weighted projective line is given by a sequence of weights and a sequence of parameters, and in this work we explain a reduction procedure using the concept of a recollement. In fact, we can classify all recollements of abelian categories which have a category of coherent sheaves of a weighted projective line in the middle.

More precisely, let $\bbX({\mathbf{p}}, \blambda)$ be a weighted projective line in the sense of Geigle--Lenzing \cite{[GL1]}, associated to a weight sequence ${{\mathbf{p}}}=(p_1, p_2,\cdots, p_t)$ and a parameter sequence $\blambda=(\lambda_1,\lambda_2,\cdots,\lambda_t)$. Here, each $p_i$ is a positive integer and $\lambda_1,\lambda_2,\cdots,\lambda_t$ are pairwise distinct points from the usual projective line. Denote by $\coh\bbX({\mathbf{p}}, \blambda)$ the category of coherent sheaves over $\bbX({\mathbf{p}}, \blambda)$.
There is a natural way to study the category $\coh\bbX({\mathbf{p}}, \blambda)$ in an inductive way, namely, by reduction/insertion of weights. The reduction of weights can be done by taking perpendicular categories \cite{[GL]}, while the insertion of weights can be achieved by the approach of $p$-cycle construction \cite{[L]} or its refinement---by expansions of abelian categories \cite{[CK]}. The above constructions can be described by the so-called reduction/insertion functors, which have nice interpretations when viewing $\coh\bbX({\mathbf{p}}, \blambda)$ as a category of certain periodic cycles.
Basing on this observation, we prove that these reduction/insertion functors form a pattern of adjoint functors, which determine a series of recollements of $\coh\bbX({\mathbf{p}}, \blambda)$. Moreover, in this way we classify the recollements of $\coh\bbX({\mathbf{p}}, \blambda)$, due to the explicit description of its localizing subcategories given by Geigle--Lenzing \cite{[GL]}.

\begin{Thm}\label{main1}(See Theorem \ref{theorem of recollement for coh})
Up to equivalences, each nontrivial recollement of $\coh\bbX({\mathbf{p}}, \blambda)$ has the following form for some sequences 
${\mathbf{q}}$ and $\underline{{\mathbf{q}}}$:
$$
\xymatrix@C=0.5cm{
\mod A(\underline{{\mathbf{q}}})\ar[rrr]^{}
&&& \coh\bbX({\mathbf{p}}, \blambda) \ar[rrr]^{\psi^{\underline{{\mathbf{q}}}}} \ar @/_1.5pc/[lll]  \ar
 @/^1.5pc/[lll]
 &&& \coh\bbX({\mathbf{p}}-{\mathbf{q}}, \blambda).
\ar @/_1.5pc/[lll]_{\psi_{\underline{{\mathbf{q}}}-\underline{{\mathbf{1}}}}} \ar
 @/^1.5pc/[lll]^{\psi_{\underline{{\mathbf{q}}}}}
 }
$$
 \end{Thm}
Here, ${{\mathbf{q}}}:=(q_1, q_2,\cdots, q_t)$ is a sequence of non-negative integers satisfying $q_i<p_i$ for each $i$ and
$\underline{{\mathbf{q}}}:=(\underline{q_1}, \underline{q_2},\cdots, \underline{q_t})$, where each $\underline{q_i}$ is a sequence of integers $(i_1, i_2, \cdots, i_{q_i})$ with $0\leq i_1<i_2<\cdots<i_{q_i}<p_i$.
In this context, we call a recollement nontrivial if its right-hand term is non-zero. The functor $\psi^{\underline{{\mathbf{q}}}}$ is called a reduction and the functor $\psi_{\underline{{\mathbf{q}}}}$ is called an insertion, each of them determines the above recollement.
Moreover, the kernel of $\psi^{\underline{{\mathbf{q}}}}$ is equivalent to a module category $\mod A(\underline{{\mathbf{q}}})$,
where $A(\underline{{\mathbf{q}}})$ is a product of path algebras of Dynkin type $\mathbb{A}$$_n$'s, and the
three remaining unlabeled functors have explicit descriptions relying on $\underline{{\mathbf{q}}}$, we refer to Theorem \ref{theorem of recollement for coh} for details.

The functors $\psi^{\underline{{\mathbf{q}}}}$ and $\psi_{\underline{{\mathbf{q}}}}$ between the categories of coherent sheaves can be lifted to their bounded derived categories, which we denote by $\Psi^{\underline{{\mathbf{q}}}}$ or $\Psi_{\underline{{\mathbf{q}}}}$ respectively.
These derived functors form an infinite sequence of adjoint pairs and yield a periodic ladder.

\begin{Thm} (See Theorem \ref{main theorem for der cat}) 
For any sequences ${\mathbf{q}}$ and $\underline{{\mathbf{q}}}$ as in Theorem \ref{main1},
there is an infinite ladder as follows, which is periodic of period ${\rm l.c.m.}(p_1,p_2,\cdots,p_t)$:
\[
\xymatrix@C=0.5cm{
D^b(\mod A(\underline{{\mathbf{q}})})
\ar@/^-3pc/[rrr]_{\vdots}%|{T_{\underline{{\mathbf{q}}}+\underline{{\mathbf{1}}}}}
\ar @/^3pc/[rrr]^{\vdots}%|{T_{\underline{{\mathbf{q}}}-\underline{{\mathbf{1}}}}}
\ar[rrr]%|{T_{\underline{{\mathbf{q}}}}}
&&& D^b(\coh\bbX({\mathbf{p}}, \blambda))
\ar@/_1.5pc/[lll]%|{T^{\underline{{\mathbf{q}}}-\underline{{\mathbf{1}}}}}
\ar@/^1.5pc/[lll]%|{T^{\underline{{\mathbf{q}}}}}
\ar[rrr]|{\Psi^{\underline{{\mathbf{q}}}}}
\ar@/^-3pc/[rrr]_{\vdots}|{\Psi^{\underline{{\mathbf{q}}}+\underline{{\mathbf{1}}}}}
\ar @/^3pc/[rrr]^{\vdots}|{\Psi^{\underline{{\mathbf{q}}}-\underline{{\mathbf{1}}}}}
&&& D^b(\coh\bbX({\mathbf{p}}-{\mathbf{q}}, \blambda)).
\ar@/_1.5pc/[lll]|{\Psi_{\underline{{\mathbf{q}}}-\underline{{\mathbf{1}}}}}
\ar@/^1.5pc/[lll]|{\Psi_{\underline{{\mathbf{q}}}}}
}
 \]
\end{Thm}

Ladders play important roles in many context.
The authors in \cite{[AKLY]} used ladders to study the derived simplicity of the derived module categories between different levels---bounded or unbounded derived categories, and finitely generated or general module categories.
Recently, a sequence of adjoint functors (viewed as a `half ladder') between compactly generated tensor-triangulated categories has been studied in \cite{[BDS]}, in order to clarify the relationship between the Grothendieck duality in algebraic geometry and the Wirthmuller isomorphism in stable homotopy theory.
In the appendix of \cite{[K2]} Keller also used an infinite sequence of adjoint functors to define the morphism enhancement of triangulated categories. We remind that a recollement of a triangulated category with Serre duality can be extended to a ladder by using the reflecting approach introduced by J{\o}gensen \cite{[J]}. This approach is useful for our further study on the stable category of vector bundles over a weighted projective line.

We denote by $\vect\bbX({\mathbf{p}}, \blambda)$ the full subcategory of $\coh\bbX({\mathbf{p}}, \blambda)$ consisting of vector bundles. In \cite{[KLM]}, Kussin, Lenzing and Meltzer proved that the subcategory $\vect\bbX({\mathbf{p}}, \blambda)$ carries a Frobenius exact structure with all the line bundles as the indecomposable projective-injectives. The associated stable category $\underline{\vect}\bbX({\mathbf{p}}, \blambda)$ is a triangulated category, which is fractional Calabi--Yau and related to singularity theory and Happel--Seidel symmetry. Moreover, in \cite{[KLM2]} they established a connection between this stable category and the problem of invariant subspaces of nilpotent linear operators, which has already been investigated by Ringel--Schmidmeier \cite{[RS]} and Simson \cite{[Sim]}. We show that the reduction/insertion functors restrict to exact functors between the stable categories of vector bundles, denoted by $\bar{\psi}^{\underline{q}}$ and $\bar{\psi}_{\underline{q}}$ respectively, which yield an periodic infinite ladder. In particular, if we are dealing with the weight triple cases, i.e, ${\mathbf{p}}=(p_1,p_2,p_3)$, we obtain a `symmetric' ladder as follows, which covers the result of \cite{[C]}.

\begin{Thm} (See Theorem \ref{main theorem for stab cat}) Let $\bbX(p_1,p_2,p_3)$ be a weighted projective line of weight type $(p_1,p_2,p_3)$. For any $0\leq q<p_3$, let $\underline{q}=(1,2,\cdots, q)$ and $\underline{q}^c=(q+1,q+2,\cdots,p_3-1)$.
Then the following diagram is an infinite ladder of period $p_3$:
$$
\xymatrix@C=0.5cm{
  \underline{\vect}\bbX(p_1,p_2,q+1)
\ar@/_3pc/[rrr]_{\vdots}|{\bar{\psi}_{\underline{q}^c+\underline{1}}}
\ar@/^3pc/[rrr]^{\vdots}|{\bar{\psi}_{\underline{q}^c-\underline{1}}}
%  \ar @/^-3.5pc/[rrr]|{\bar{\psi}_{({\underline{q}}-i_1{\underline{1}})^c}} \ar @/^3.5pc/[rrr]_{}|{\bar{\psi}_{({\underline{q}}-i_1{\underline{1}})^c}}
\ar[rrr]|{\bar{\psi}_{{\underline{q}}^c}}
&&& \underline{\vect}\bbX(p_1,p_2,p_3)
\ar@/_1.5pc/[lll]|{\bar{\psi}^{{\underline{q}}^c-{\underline{1}}} }
\ar@/^1.5pc/[lll]|{\bar{\psi}^{q^c}}
\ar[rrr]|{\bar{\psi}^{\underline{q}}}
%\ar@/^-3.5pc/[rrr]|{\bar{\psi}^{{\underline{q}}+(p-i_1){\underline{1}}}=\bar{\psi}^{{\underline{q}}-i_1{\underline{1}}}}
%\ar @/^3.5pc/[rrr]_{}|{\bar{\psi}^{{\underline{q}}-i_1{\underline{1}}}}
\ar@/_3pc/[rrr]_{\vdots}|{\bar{\psi}^{{\underline{q}}+{\underline{1}}}}
\ar@/^3pc/[rrr]^{\vdots}|{\bar{\psi}^{{\underline{q}}-{\underline{1}}}}
&&& \underline{\vect}\bbX(p_1,p_2,p_3-q).
\ar@/_1.5pc/[lll]|{\bar{\psi}_{{\underline{q}}-{\underline{1}}}}
\ar@/^1.5pc/[lll]|{\bar{\psi}_{\underline{q}}}
}
$$
\end{Thm}

This ladder above has direct applications in tilting theory for the stable category of vector bundles over a weighted projective line. In fact, we can prove \cite[Theorem B]{[KLM]} in an inductive way. This also provides an effective method to construct new tilting objects in the stable category $\underline{\vect}\bbX({\mathbf{p}}, \blambda)$, due to the explicit expressions of the functors in the ladder above. These will be investigated in the forthcoming paper.

This paper is organized as follows. Some preliminaries about recollements and ladders are collected in Section 2. We also recall the $p$-cycle construction for the category of coherent sheaves over an exceptional curve there.
In Section 3 we define the reduction functor $\psi^i$ and the insertion functor $\psi_i$ for any integer $i$, and investigate their basic properties. We prove that these functors form an infinite sequence of adjoint functors, which can be used to construct recollements of the category of coherent sheaves and ladders of its bounded derived category. The applications to weighted projective lines are discussed in Section 4. We  classify all the recollements for the category of coherent sheaves over a weighted projective line, which induces a periodic infinite ladder of its bounded derived category. Finally, we show that the reduction/insertion functors restrict to exact functors between the stable categories of vector bundles, yielding an infinite ladder which is also periodic.

\section{Preliminary}

In this section, we fix our notations and collect some basic facts about recollements and ladders, and introduce the category of coherent sheaves on a noncommutative exceptional curve by using $p$-cycles construction due to Lenzing \cite{[L]}.

Throughout this paper we fix an arbitrary field $\mathbf{k}$, and denote by $D:=\Hom_{\mathbf{k}}(-, \mathbf{k})$. For an algebra $A$ we denote by $\mod A$ the category of finitely generated right $A$ modules. Let $F: \cC\to \cC'$ be an additive functor between additive categories. The \emph{kernel} $\Ker(F)$ of $F$ is the full subcategory of $\cC$ formed by all objects $C$ such that $FC=0$. The \emph{essential image} $\im(F)$ is the full subcategory of $\cC'$ formed by all objects $C'$ such that $C'$ is isomorphic to $FC$ for some $C\in\cC$. The functor $F: \cC\to\cC'$ is called a \emph{quotient functor} if
$F$ induces an equivalence $\cC[\Sigma_F^{-1}]\cong \cC'$, where $\Sigma_F$ denotes all the classes of morphisms $f$ in $\cC$ with $F(f)$ invertible, and $\cC[\Sigma_F^{-1}]$ is the localization category of $\cC$ with respect to $\Sigma_F$. A sequence of additive categories and additive functors
 $\xymatrix@C=1cm{
\cc'\ar[r]^{F_1}
& \cc
\ar[r]^{F_2}
& \cc''
}$ is called \emph{exact} if the following conditions hold:
\begin{itemize}
\item[(1)] $F_1$ is fully faithful;
\item[(2)] $F_2$ is a quotient functor;
\item[(3)] $\Ker(F_2)=\im(F_1)$.
\end{itemize}

\subsection{Adjoint pair}

Let $\cC$ and $\cC'$ be arbitrary categories, and let
\begin{eqnarray}\label{adjoint pair diagram} \xymatrix@!=3pc{\cC
\ar@<1ex>[r]^{F} &\cC'
\ar@<1ex>[l]^{G}
}
\end{eqnarray}
be a pair of functors between them. Then $(F,G)$ is called an \emph{adjoint pair}
if there exists a natural isomorphism
$$\Hom_{\cC'}(FC, C')\cong\Hom_{\cC}(C, GC'),$$ which are functorial in objects $C\in\cC$ and $C'\in\cC'$. An adjoint pair $(F,G)$ is often denoted by $F\dashv G$, and call $F$ a left adjoint and $G$ a right adjoint. The following three statements are equivalent for an adjoint pair $(F,G)$ (see \cite[Prop. I.1.3]{[GZ]}):
\begin{itemize}
 \item [(1)] $FG\cong \id_{\cC'}$;
  \item [(2)] $G$ is fully faithful;
   \item [(3)] $F$ is a quotient functor.
\end{itemize}

It is well known that adjoint pairs are compatible with auto-equivalences of categories and also with compositions. More precisely, if $(F,G)$ is an adjoint pair between $\cC$ and $\cC'$, and there exist auto-equivalences $\sigma$ of $\cC$ and $\sigma'$ of $\cC'$ respectively, then $(\sigma' F\sigma, \sigma^{-1}G (\sigma')^{-1})$ is also an adjoint pair; if $(F_1, G_1)$ and $(F_2, G_2)$ are adjoint pairs as follows
\begin{eqnarray}\xymatrix@!=3pc{\cC
\ar@<1ex>[r]^{F_1} &\cC'\ar@<+1ex>[r]^{F_2}
\ar@<1ex>[l]^{G_1}  &
\cC'',\ar@<1ex>[l]^{G_2}
}\label{eq:recollement}
\end{eqnarray}
then $(F_2F_1, G_1G_2)$ is also an adjoint pair.

\subsection{Recollement and Ladder}
Let $\cC$ be an additive category.
A {\em recollement} of $\cC$ is the following diagram
\begin{equation}\label{recollement}
\xymatrix@C=1cm{
\cc'\ar[rrr]|{i_*=i_!}
&&& \cc
\ar@/_1.5pc/[lll]|{i^*}
\ar@/^1.5pc/[lll]|{i^!}
\ar[rrr]|{j^!=j^*}
&&& \cc''
\ar@/_1.5pc/[lll]|{j_!}
\ar@/^1.5pc/[lll]|{j_*}
}
 \end{equation}
of additive categories and additive functors such that
\begin{enumerate}
\item $(i^\ast,i_\ast)$,\,$(i_!,i^!)$,\,$(j_!,j^!)$,\,$(j^\ast,j_\ast)$
are adjoint pairs;
\item  $i_\ast,\,j_\ast,\,j_!$  are fully faithful;
\item  $\im (i_\ast)=\Ker (j^{\ast})$.
\end{enumerate}
By definition, the middle row of a recollement
is an exact sequence.
Two recollements of $\cC$ are called \emph{equivalent} if the related subcategories $\Ker(j^{\ast})$ coincide.

The recollements of abelian categories or triangulated categories have been studied widely, see for examples \cite{[BBD], [P]}. Here we only recall when the recollement (\ref{recollement}) is determined by its left (right) half part.
If all of $\cC', \cC, \cC''$ are abelian categories,
then the right half recollement of (\ref{recollement}), i.e, the adjoint triple $(j_!,j^\ast,j_\ast)$ together with $j_!$ and (=or) $j_\ast$ fully faithful, can be completed to the whole recollement; however, the dual statement is not true, i.e, the left half recollement can not be completed to a recollement in general, see \cite[Rem. 2.5]{[P]}. But if we are dealing with triangulated categories, then either left or right half recollement can be completed to a recollement, see \cite[Thm. 1.1]{[CPS1]} and \cite[Thm. 2.1]{[CPS2]} respectively, which relies on the basic fact that
the left (resp. right) adjoint functors of exact functors in triangulated categories are automatically exact \cite[Lem. 5.3.6]{[N]}.

A \emph{ladder} of an additive category $\cC$ is a finite or infinite diagram of additive categories and additive functors
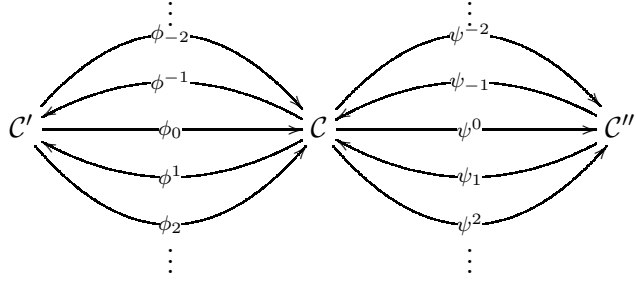
\[
\xymatrix@C=1cm{
\cc' \ar @/^-3pc/[rrr]|{\phi_{2}}_{\vdots}
\ar @/^3pc/[rrr]_{}|{\phi_{-2}} ^{\vdots}
\ar[rrr]|{\phi_0}
&&& \cc
\ar@/_1.5pc/[lll]|{\phi^{-1}}
\ar@/^1.5pc/[lll]|{\phi^{1}}
\ar[rrr]|{\psi^{0}}
\ar @/^-3pc/[rrr]|{\psi^{2}}_{\vdots}
\ar @/^3pc/[rrr]_{}|{\psi^{-2}}^{\vdots}
&&& \cc''
\ar@/_1.5pc/[lll]|{\psi_{-1}}
\ar@/^1.5pc/[lll]|{\psi_{1}}
}
 \]
such that any three consecutive rows form a recollement. The rows are labelled by a subset of $\bbZ$ and multiple occurrence of the same recollement is allowed. The \emph{height} of a ladder is the number of recollements contained in it (counted with multiplicities). A recollement is considered to be a ladder of height 1.
A ladder is \emph{periodic} if there exists a positive integer $n$ such that the $n$-th recollement going upwards (respectively, going downwards) in the ladder is equivalent to the recollement which is considered to be a ladder of height one. The minimal such positive integer $n$ is the \emph{period} of the ladder.
If the middle term of a recollement of triangulated categories has Serre duality, then by using the reflecting approach induced by J{\o}gensen \cite{[J]} we can obtain a ladder.

\subsection{The category of $p$-cycles}

Let $\cal{H}$ be a category of coherent sheave on a noncommutative exceptional curve $\bbX$. Recall that the torsion subcategory $\cH_0$ of $\cal{H}$ has a decomposition $\cH_0=\coprod_{x\in \bbX}\cU_x$, where $\cU_x$ is consisting of coherent sheaves concentrated in $x$. Let $\cS_x$ be the semisimple subcategory of $\cU_x$. Each point $x\in\bbX$ determines (by means of a mutation with respect to $\cS_x$) a functor $\sigma_x: \cH\to \cH; \; E\mapsto E(x)$, together with a natural transformation $x: Id\to \sigma_x$, also denoted by the symbol $x$. For each sheaf $E$ with $\Hom(\cU_x, E)=0$, these data are given by the $\cS_x$-universal extension:
\begin{equation}\label{universal extension}\xymatrix{0\ar[r]& E\ar[r]^{x_E}& E(x)\ar[r]& E_x \ar[r]& 0,}
\end{equation}
where
$E_x$ belongs to $\cS_x$. In particular, if $E$ is a torsion sheaf concentrated in a point different from $x$, then $E(x)=E$ and $x_E=\id_E$. Moreover, if $E$ is a torsion sheaf concentrated in $x$ then $E(x)=\tau^{-1}E$, where $\tau$ is the Auslander-Reiten translation of $\cal{H}$, and the kernel of $x_E$ equals the simple socle of $E$.

Fix a point $x$ of $\bbX$ and an integer $p\geq 1$. A \emph{$p$-cycle} $E$ concentrated in $x$ is a diagram
$$\xymatrix{
\cdots\ar[r] & E_{n}\ar^{x_{n}}[r] &  E_{n+1}\ar^{x_{n+1}}[r] &\cdots \ar[r] & E_{n+p}\ar^{x_{n+p}}[r]&\cdots},$$
which is $p$-periodic in the sense that $E_{n+p}=E_{n}(x)$, $x_{n+p}=x_{n}(x)$ and moreover $x_{n+p-1}\cdots x_{n+1} x_{n}=x_{E_{n}}$ holds for each integer $n$. A morphism $u: E\to F$ of $p$-cycles concentrated in the same point $x$ is a sequence of morphism $u_n: E_n\to F_n$ which is $p$-periodic, i.e, $u_{n+p}=u_{n}(x)$ for each $n$ and such that each diagram
$$\xymatrix{
E_n\ar[r]^{x_n}\ar[d]^{u_n} &E_{n+1}\ar[d]^{u_{n+1}}\\
F_n\ar[r]^{x_n}& F_{n+1}
}$$
commutes.
We denote $p$-cycles in the form
\begin{equation}\label{p-cycle}\xymatrix{
E_0\ar[r]^{x_0} & E_{1}\ar[r]^{x_{1}} & \cdots \ar[r]^{x_{p-2}}&  E_{p-1}\ar[r]^{x_{p-1}} & E_{0}(x),
}\end{equation}
and denote the category of all $p$-cycles concentrated in $x$ by $\cH(p;x)$. %{\red which will be simplified 

According to \cite{[L]}, the category $\cH(p;x)$ is connected, abelian and noetherian, where exactness and formation of kernels and cokernels have a pointwise interpretation. Moreover, it is equivalent to the category of coherent sheaves on a smooth projective curve.
Note that any functor $\sigma$ on $\cH$ induces a functor on $\cH(p;x)$ via pointwise actions, which will be still denoted by $\sigma$.
Moreover, $\cH(p;x)$ is again equipped with a natural shift automorphism $\overline{\sigma}_{x}$, which satisfies $\overline{\sigma}_{x}^p=\sigma_x$ and sends
a $p$-cycle of the form (\ref{p-cycle}) to $$\xymatrix{
 E_{1}\ar[r]^{x_{1}} & E_2\ar[r]^{x_2} &\cdots \ar[r]^{x_{p-1}}&  E_{0}(x)
\ar^{x_{0}}[r]&  E_{1}(x).}$$

\subsection{Simple sheaves and line bundles in $\cH(p;x)$}\label{section of simples}

Recall from \cite{[L]} that there is a full exact embedding $$\iota: {\cal{H}}\to \cH(p;x); \quad E\mapsto (E=\cdots =E\to E(x)).$$
We will identify $\cal{H}$ with the resulting exact subcategory of $\cH(p;x)$.
It is important to remind that the embedding functor $\iota$ is rank-preserved.

The simple objects in $\cH(p;x)$ occur in two types:

\begin{itemize}
  \item [(1)] the simple objects of $\cal{H}$ which are concentrated in a point $y$ different from $x$;
  \item [(2)] for each simple object $S$ from $\cal{H}$, concentrated in $x$, the $p$ simples:
  \begin{itemize}
  \item[] \qquad $S_1:\; 0\to 0\to \cdots \to 0\to S\to 0;$
  \item[] \qquad $S_2:\; 0\to 0\to \cdots \to S \to 0\to 0;$
  \item[] \qquad $\qquad\qquad\qquad\cdots\cdots$
  \item[] \qquad $S_{p-1}:\; 0\to S\to \cdots \to 0 \to 0\to 0;$
  \item[] \qquad $S_{p}: \;S\to 0\to \cdots \to 0 \to 0\to S(x).$
\end{itemize}
\end{itemize}
Moreover, each $S_{i}$ is exceptional and $\End_{\cH(p;x)}(S_i)\cong\End_{\cH}(S):=\Delta$, which is a finite skew field extension of $\mathbf{k}$. Recall that $S(x)=\tau^{-1}S$. The non-zero extensions between the simples of type (2) are given by $$\Ext_{\cH(p;x)}^1(S_{i+1},S_i)\cong \Delta\cong \Ext_{\cH(p;x)}^1(S(x)_{1},S_p)$$ for $1\leq i\leq p-1$.
Consequently, $\bar{\tau}^{-1}(S_{i})=S_{i+1}$ for $1\leq i\leq p-1$ and $\bar{\tau}^{-1}(S_{p})=S(x)_1$, where $\bar{\tau}$ is the Auslander--Reiten translation of $\cH(p;x)$.
In particular, if there is a unique simple sheaf $S$ in $\cal{H}$ concentrated in $x$, i.e, $S(x)=S$, then $\bar{\tau}(S_{i+1})=S_i$, where the indices are taken modulo $p$. That is, the extension-closed subcategory generated by the simples of type (2) is a tube of rank $p$.

Let $E$ be a $p$-cycle in $\cH(p;x)$ of the form (\ref{p-cycle}). Then $E$ is a torsion sheaf (resp. bundle) if and only if each $E_i$ a torsion sheaf (resp. bundle) in $\cal{H}$. In the bundle case, each $x_i$ is a monomorphism since each $x_{E_i}: E_i\to E_i(x)$ is a monomorphism by (\ref{universal extension}).
Now we consider the line bundles in $\cH(p;x)$. Obviously, for any line bundle $L\in\cH$ and $0\leq i\leq p-1$, the $p$-cycle of the following form is a line bundles in $\cH(p;x)$:
$$\xymatrix@C=0.5cm{
\bar{\sigma}_x^{i}(\iota L):& L\ar@{=}[r] & \cdots \ar@{=}[r] &  L\ar[rr]^{x_{p-i-1}} && L(x)\ar@{=}[r] & L(x)\ar@{=}[r] &\cdots \ar@{=}[r] &L(x).}$$
In fact, any line bundle in $\cH(p;x)$ can be obtained in this way.

\begin{Prop}\label{extensions and homs}
Each line bundle in $\cH(p;x)$ has the form $\bar{\sigma}_x^{i}(\iota L)$
for some line bundle $L$ in $\cH$ and $0\leq i\leq p-1$. Moreover, if $\Ext_{\cal{H}}^1(S, L)\neq 0$, then
\begin{itemize}
  \item[(1)] $\Ext^1(S_k, \bar{\sigma}_x^{i}(\iota L))\neq 0$ if and only if $k=i+1$;
  \item[(2)] $\Hom(\bar{\sigma}_x^{i}(\iota L), S_k)\neq 0$ if and only if $k=i$.
\end{itemize}
\end{Prop}

\begin{pf}

For each line bundle $L\in\cH$, there exists a unique simple sheaf $S_L$ concentrated in $x$ satisfying $\Ext^1(S_L, L)\neq 0$. Hence the $\cS_x$-universal extension of $L$ has the form:
$$\xymatrix{0\ar[r]& L\ar[r]^{x_L}& L(x)\ar[r]& S_L \ar[r]& 0.}$$

Now let $\xymatrix{
L_0\ar[r]^{x_0} & L_{1}\ar[r]^{x_{1}} & \cdots \ar[r]^{x_{p-2}}& L_{p-1}\ar[r]^{x_{p-1}} & L_{0}(x)}$ be a line bundle in $\cH(p;x)$, then each $L_i$ is a line bundle in $\cH$ and each $x_i$ is injective.
We claim that all but one of these $x_i's$ are isomorphisms. In fact, assume $x_k$ is not an isomorphism for some $k$. We set $f=x_{k-1}\cdots x_1x_0: L_0\to L_{k}$ and $g=x_{p-1}\cdots x_{k+1}x_k: L_k\to L_{0}(x)$, then $g$ is not an isomorphism. Consider the following commutative diagram:
$$\xymatrix@C=0.5cm{
 0\ar[r]& L_0\ar[r]^{f} \ar@{=}[d]& L_k \ar[r] \ar^{g}[d]& S' \ar[r]\ar[d]^{h} &0\\
0\ar[r]& L_0\ar[r]^{x_{L_0}} & L_0(x) \ar[r]& S_{L_0} \ar[r]&0.}$$
Then $g$ is not an isomorphism implies so as $h$, which means $S'=0$ since $S_{L_0}$ is simple. It follows that $f$ is an isomorphism, and so as $x_i$ for $0\leq i\leq k-1$. Similarly, we can prove $x_i$ is an isomorphism for $k+1\leq i\leq p-1$. This proves the first part of the result.

The second statement follows from the definition of morphisms between $p$-cycles, which are involved in the following diagram:
$$\xymatrix@C=0.5cm{
\bar{\sigma}_x^{k-1}(\iota L):\ar[d]& L\ar@{=}[r] \ar[d]& \cdots \ar@{=}[r] \ar[d]& L\ar@{=}[r] \ar[d]& L\ar[r]^{x_{p-k}} \ar[d]^{x_{L}}& L(x)\ar@{=}[r] \ar[d]&\cdots \ar@{=}[r]\ar[d] &L(x)\ar[d]\\
\bar{\sigma}_x^{k}(\iota L):\ar[d]& L\ar@{=}[r] \ar[d]& \cdots \ar@{=}[r]\ar[d] &  L\ar[r]^{x_{p-k-1}} \ar[d]& L(x)\ar@{=}[r]\ar[d] & L(x)\ar@{=}[r] \ar[d]&\cdots \ar@{=}[r]\ar[d] &L(x)\ar[d]\\
S_k:& 0\ar@{=}[r] & \cdots \ar@{=}[r] &  0\ar[r] & S\ar[r] & 0\ar@{=}[r] &\cdots \ar@{=}[r] &0.}$$

\end{pf}

Let ${\cal{S}}=\{S_2,S_3,\cdots, S_{p}\,|\,$ S simple in $\cal{H}$ and concentrated in $x \}$, then the extension closure $\langle \cal{S}\rangle$ of $\cS$ is localizing in $\cH(p;x)$,
and the category $\langle\cal{S}\rangle^{\perp}$ consists of exactly those $p$-cycles of the form (\ref{p-cycle}) such that $x_0, x_1,\cdots, x_{p-2}$ are isomorphisms, hence---up to isomorphism---agree with the objects from $\cal{H}$. It follows that
$$\cH(p;x)/\langle \cal{S}\rangle\cong \langle\cal{S}\rangle^{\perp}\cong \cal{H}.$$

\section{Recollements and ladders induced by $p$-cycle construction}

In this section, we will investigate the recollements and ladders induced by $p$-cycle construction for noncommutative exceptional curves. Throughout this section, let $\bbX$ be a noncommutative exceptional curve, and let $\cal{H}$ be the category of coherent sheave over $\bbX$. We fix a point $x$ of $\bbX$. For any integer $p\geq 1$, denote by $\cH(p):=\cH(p;x)$ the category of all the $p$-cycles in $\cH$ which are concentrated in $x$.

\subsection{Reduction functor $\psi^i$ and insertion functor $\psi_i$}

For any $0\leq i\leq p-1$, we define a pair of functors $(\psi^i, \psi_i)$ as follows.
The \emph{reduction functor} $\psi^i: \cH(p)\to \cH(p-1)$ is defined by sending a $p$-cycle $E$ of the form (\ref{p-cycle})
to a $(p-1)$-cycle:
$$\xymatrix{
E_0\ar[r]^{x_0} & \cdots\ar[r] & E_{i-1}\ar[r]^{x_ix_{i-1}} &  E_{i+1}\ar[r]^{x_{i+1}} &\cdots\ar[r] & E_{p-1}\ar[r]^{x_{p-1}} & E_{0}(x),
}
$$
and sending a morphism of $p$-cycles $(u_0,u_1,\cdots, u_{p-1})$ to
$(u_0,\cdots, u_{i-1}, u_{i+1}, \cdots, u_{p-1}).$
Conversely, the \emph{insertion functor} $\psi_i: \cH(p-1)\to \cH(p)$ is defined by sending a $(p-1)$-cycle $E'$ of the form
\begin{equation}\label{p-1 cycle}
\xymatrix{
E'_0\ar[r]^{x'_0} & E'_{1}\ar[r]^{x'_{1}} &\cdots\ar[r] & E'_{p-2}\ar[r]^{x'_{p-2}} & E'_{0}(x)
}\end{equation}
to a $p$-cycle
\begin{equation}\label{image of p-1 cycle}\xymatrix@C=0.6cm{
E'_0\ar[r]^{x'_0} & \cdots\ar[r] & E'_{i}\ar@{=}[r]  &E'_{i}\ar[r]^{x'_i}  & E'_{i+1}\ar[r]^{x'_{i+1}} &\cdots\ar[r] & E'_{p-2}\ar[r]^{x'_{p-2}} & E'_{0}(x),}
 \end{equation}
 and sending a morphism of $(p-1)$-cycles $(u'_0,\cdots, u'_{i}, u'_{i+1}, \cdots, u'_{p-2})$ to a morphism of $p$-cycles $(u'_0,\cdots, u'_{i}, u'_{i}, u'_{i+1}, \cdots, u'_{p-2})$.

Roughly speaking, $\psi^i(E)$ is obtained from $E$ by deleting the item $E_i$ and taking the composition of the morphisms $\xymatrix{
E_{i-1}\ar[r]^{x_{i-1}}&E_{i}\ar[r]^{x_{i}} & E_{i+1}}$; conversely, $\psi_i(E')$ is obtained from $E'$ by inserting in the $i$-th position a copy of $E'_i$: $\xymatrix{
E'_{i}\ar@{=}[r]  &E'_{i}\ar[r]^{x'_i}  & E'_{i+1}}$.

Recall that the functor $\sigma_x: \cH\to \cH; \;E\mapsto E(x)$ induces a functor on $\cH(p)$ (resp. $\cH(p-1)$) via pointwise actions, which shares the same notation $\sigma_x$. Denote by $\overline{\sigma}_{x}$ and $\overline{\sigma}_{x}'$ the natural shift automorphisms on $\cH(p)$ and $\cH(p-1)$ respectively, then $$\overline{\sigma}_{x}^p=\sigma_x: \cH(p)\to \cH(p)\quad{\text {and}}\quad(\overline{\sigma}_{x}')^{p-1}=\sigma_x: \cH(p-1)\to\cH(p-1).$$
By the above constructions, we have $$\psi^{p-1}=\psi^0\overline{\sigma}_{x}^{-1}
\quad{\text {and}}\quad \psi_{p-1}=\overline{\sigma}_{x}\psi_0.$$
For any $n\in\bbZ$ and $0\leq i\leq p-1$, define
\begin{equation}\label{def of fun for gen case} \psi^{np+i}=(\overline{\sigma}_{x}')^{-n}\psi^{i}
\text{\quad and \quad}\psi_{np+i}=\psi_{i}(\overline{\sigma}_{x}')^{n}.
\end{equation}

\begin{Prop}\label{conjugate} For any $i\in\bbZ$, $\psi^{i}=(\overline{\sigma}_{x}')^{-1}\psi^{i-1}\overline{\sigma}_{x}$ and
$\psi_{i}=\overline{\sigma}_{x}^{-1}\psi_{i-1}\overline{\sigma}_{x}'$.
\end{Prop}

\begin{pf}
We only prove the first formula. The second one can be proved similarly.
For $1\leq i\leq  p-1$, it follows immediately from the construction that $\psi^{i}=(\overline{\sigma}_{x}')^{-1}\psi^{i-1}\overline{\sigma}_{x}$. Moreover, by (\ref{def of fun for gen case}) we have
$\psi^{-1}=\overline{\sigma}_{x}'\psi^{p-1}=\overline{\sigma}_{x}'\psi^0\overline{\sigma}_{x}^{-1}$, hence $\psi^{0}=(\overline{\sigma}_{x}')^{-1}\psi^{-1}\overline{\sigma}_{x}$.
In general case, we have $$\psi^{np+i}=(\overline{\sigma}_{x}')^{-n}\psi^{i}
=(\overline{\sigma}_{x}')^{-n}(\overline{\sigma}_{x}')^{-1}\psi^{i-1}\overline{\sigma}_{x}
=(\overline{\sigma}_{x}')^{-1}\psi^{np+i-1}\overline{\sigma}_{x}.$$
\end{pf}

\subsection{Adjointness between reduction and insertion functors}

Now we state our main observation on the adjointness between the functors $\psi^i$'s and $\psi_i$'s, which plays a key role in this paper.

 \begin{Prop}\label{adjoint sequence in abelian category}
For any $i\in\bbZ$, $(\psi_{i-1}, \psi^i, \psi_i)$ form an adjoint triple. Consequently, both $\psi^i$ and $\psi_i$ are exact functors and there is an infinite sequence of adjoint functors
\begin{equation}\label{adjoint sequence of psi} \cdots\dashv\psi^{-1}\dashv\psi_{-1}\dashv\psi^0\dashv\psi_0\dashv\psi^1\dashv \cdots  .
\end{equation}
\end{Prop}

\begin{pf}
We first show that $\psi^i$ admits a right adjoint functor $\psi_i$. By (\ref{def of fun for gen case}), we can restrict to $0\leq i\leq p-1$.
Let $E$ be a $p$-cycle of the form (\ref{p-cycle}) and $E'$ be a $(p-1)$-cycle of the form (\ref{p-1 cycle}).
Then $\Hom(\psi^i(E), E')$ consists of sequences of morphisms $(u_0, \cdots, u_{i-1}, u_{i+1}, \cdots, u_{p-1})$ satisfying that the following diagram commutes:
$$\xymatrix{
E_0\ar[r]^{x_0}\ar[d]^{u_0} & \cdots\ar[r] & E_{i-1}\ar[r]^{x_ix_{i-1}}\ar[d]^{u_{i-1}} &  E_{i+1}\ar[r]^{x_{i+1}} \ar[d]^{u_{i+1}} &\cdots\ar[r] & E_{p-1}\ar[r]^{x_{p-1}}\ar[d]^{u_{p-1}} & E_{0}(x)\ar[d]^{u_{0}(x)}
\\
E'_0\ar[r]^{x'_0} & \cdots\ar[r] & E'_{i-1}\ar[r]^{x'_{i-1}} &  E'_{i}\ar[r]^{x'_{i}} &\cdots\ar[r] & E'_{p-2}\ar[r]^{x'_{p-2}} & E'_{0}(x).
}$$
Similarly, $\Hom(E, \psi_i(E'))$ consists of sequences $(u'_0, \cdots, u'_{i-1}, u'_i, u'_{i+1}, \cdots, u'_{n})$ such that the following diagram commutes:
$$\xymatrix{
E_0\ar[r]^{x_0}\ar[d]^{u'_0} & \cdots\ar[r] &E_{i-1}\ar[r]^{x_{i-1}}\ar[d]^{u'_{i-1}}& E_{i}\ar[r]^{x_i}\ar[d]^{u'_{i}} &  E_{i+1}\ar[r]^{x_{i+1}} \ar[d]^{u'_{i+1}} &\cdots\ar[r] & E_{p-1}\ar[r]^{x_{p-1}}\ar[d]^{u'_{p-1}} & E_{0}(x)\ar[d]^{u'_{0}(x)}
\\
E'_0\ar[r]^{x'_0} & \cdots\ar[r] & E'_{i-1}\ar[r]^{x'_{i-1}} &  E'_{i}\ar@{=}[r] &E'_{i}\ar[r]^{x'_{i}} &\cdots\ar[r] & E'_{p-2}\ar[r]^{x'_{p-2}} & E'_{0}(x),
}$$
which then implies $u'_i=u'_{i+1}x_i$ (i.e, $u'_i$ is determined by $u'_{i+1}$).
Therefore, the map $$(u_0, \cdots, u_{i-1}, u_{i+1}, \cdots, u_{p-1})\mapsto(u_0, \cdots, u_{i-1}, u_{i+1}x_i, u_{i+1}, \cdots, u_{p-1})$$ induces a natural isomorphism $$\Hom_{\cH(p-1)}(\psi^i(E), E')\cong\Hom_{\cH(p)}(E, \psi_i(E')),$$ which is functorial in both $E$ and $E'$ by standard arguments. This proves the adjointness of $(\psi^i, \psi_i)$.

By similar arguments we can show the adjointness $\psi_{i-1}\dashv\psi^i$.
It follows that $\psi^i$ (resp. $\psi_i$) is exact since it admits both a left and a right adjoint.
\end{pf}

\begin{Lem}\label{exact of r and l} For any $i\in\bbZ$, $\psi^i\psi_i=\id_{\cH(p-1)}=\psi^i\psi_{i-1}$. In particular, each $\psi_i$ is fully faithful, and each $\psi^i$ is a quotient functor.
\end{Lem}

\begin{pf} For any $0\leq i\leq p-1$, by construction we have $\psi^i\psi_i=\id_{\cH(p-1)}$. It follows that $$\psi^{np+i}\psi_{np+i}
=(\overline{\sigma}_{x}')^{-n}\psi^{i}\psi_{i}(\overline{\sigma}_{x}')^{n}
=(\overline{\sigma}_{x}')^{-n}\id_{\cH(p-1)}(\overline{\sigma}_{x}')^{n}=\id_{\cH(p-1)}.$$
Similarly, we can prove $\psi^i\psi_{i-1}=\id_{\cH(p-1)}$ for any $i\in\bbZ$.
It follows that $\psi_i$ is fully faithful and $\psi^i$ is a quotient functor, see for example \cite[Prop. I.1.3]{[GZ]}.
\end{pf}

Recall that the set of simple sheaves in $\cal{H}$ concentrated in $x$ is denoted by $\cS_x$. For any $i\in\bbZ$, set $\tilde{\mathbf{S}}_{i}=\bigoplus_{S\in \cS_x}S_i$, where the indices are taken modulo $p$.
\begin{Lem}\label{ker and image} For any $i\in\bbZ$, $\Ker(\psi^i)=\add(\tilde{\mathbf{S}}_{p-i})$ and $\im(\psi_i)=\tilde{\mathbf{S}}_{p-i}^{\perp}$.
\end{Lem}

\begin{pf}
By definition we have $\Ker(\psi^{np+i})=\Ker(\psi^i)$ and $\im(\psi_{np+i})=\im(\psi_i)$. Hence we can assume $0\leq i\leq p-1$. By the construction of the functor $\psi^i$, it is obvious that $\Ker(\psi^i)$ consists of those $p$-cycles of the form (\ref{p-cycle}) with $E_j=0$ for any $j\neq i$. It follows that $\Ker(\psi^i)=\add(\tilde{\mathbf{S}}_{p-i})$ by the explicit description of simples in $\cH(p)$ in Section \ref{section of simples}.

In the following we remain to show that $\im(\psi_i)=\tilde{\mathbf{S}}_{p-i}^{\perp}$ for $0\leq i\leq p-1$. 
Firstly, we show $\im(\psi_i)\subseteq \tilde{\mathbf{S}}_{p-i}^{\perp}$. Assume $E\in \im(\psi_i)$, then $E$ has the form (\ref{image of p-1 cycle}).
For any simple sheaf $S$ in $\cal{H}$ concentrated in $x$, recall that $S_{p-i}$ has the form $ 0\to \cdots \to S\to \cdots \to 0$, where $S$ sits in the position $i$. Hence $\Hom(S_{p-i}, E)$ is determined by a morphism $f_i\in\Hom(S, E'_{i})$ with the following commutative diagram,
$$\xymatrix{
S\ar[r] \ar[d]^{f_i}&  0 \ar[d]\\ E'_{i}\ar@{=}[r] &  E'_{i}.
}$$
It follows that $f_i=0$ and then $\Hom(S_{p-i}, E)=0$.
Similarly, we can prove that $\Hom(E, S_{p-i-1})=0$. Then by Serre duality, $\Ext^1(S_{p-i}, E)=D\Hom(E, S_{p-i-1})=0$. It follows that $E\in S_{p-i}^{\perp}$ and then $\im(\psi_i)\subseteq \tilde{\mathbf{S}}_{p-i}^{\perp}$.

On the other hand, for any indecomposable $E\in \tilde{\mathbf{S}}_{p-i}^{\perp}$, assume $E$ has the form (\ref{p-cycle}).
We are going to prove that $x_i$ is an isomorphism. If $E$ is a torsion sheaf concentrated in a point different from $x$, then all of $x_0, x_1,\cdots, x_{p-1}$ are isomorphisms, we are done. If else, we consider the following exact sequence in $\cH$: $$\xymatrix{
0\ar[r]& \Ker(x_i)\ar[r]^-{\vez}&  E_i\ar[r]^-{x_i} &  E_{i+1} \ar[r]^-{\pi}&\Coker(x_i)\ar[r]&0.
}$$
We claim both of $\Ker(x_i)$ and $\Coker(x_i)$ are torsion sheaves concentrated in $x$. In fact, if $E$ is a torsion sheaf in $\cH(p)$ concentrated in $x$, then both of $E_{i}$ and $E_{i+1}$ are torsion sheaves in $\cH$ concentrated in $x$, hence so as $\Ker(x_i)$ and $\Coker(x_i)$; if $E$ is a vector bundle, then by \cite[Sect. 4.2]{[L]} the $p$-cycle (\ref{p-cycle}) can be interpreted as a filtration $$E_0\subseteq E_1\cdots\subseteq E_{p-1}\subseteq E_0(x),$$ or equivalently as a filtration
$$0=E_0/E_0\subseteq E_1/E_0\cdots\subseteq E_{p-1}/E_0\subseteq E_0(x)/E_0$$ of the fibre $E_x=E_0(x)/E_0$ of $E$ at $x$ given by the universal extension (\ref{universal extension}). It follows that $\Ker(x_i)=0$ and $\Coker(x_i)$ is a factor of $E_x$, hence concentrated in $x$. This proves the claim.

Now if $x_i$ is not injective, then $\vez\neq 0$, which gives to a nonzero morphism in $\Hom(\tilde{\mathbf{S}}_{p-i}, E)$ via the following commutative diagram on the left-hand side, contradicting to $E\in \tilde{\mathbf{S}}_{p-i}^{\perp}$.
$$\xymatrix{
\Ker(x_i)\ar[r] \ar[d]_{\iota}&  0 \ar[d]
&&&& E_i\ar[r]^{x_i} \ar[d]&  E_{i+1} \ar[d]^{\pi}
\\ E_{i}\ar[r]_{x_i} &  E_{i+1}
&&&& 0 \ar[r] &  \Coker(x_i) .
}$$
Similarly, if $x_i$ is not surjective, then $\pi\neq 0$, which gives to a nonzero morphism in $\Hom(E, \tilde{\mathbf{S}}_{p-i-1})\cong D\Ext^{1}(\tilde{\mathbf{S}}_{p-i}, E)$ via the above commutative diagram on the right-hand side, another contradiction. Therefore, $x_i$ is an isomorphism, which means that $E\in \im (\psi_i)$ and then $\tilde{\mathbf{S}}_{p-i}^{\perp}\subseteq\im (\psi_i)$. Therefore, $\im (\psi_i)=\tilde{\mathbf{S}}_{p-i}^{\perp}$, we are done.
\end{pf}

\subsection{Recollements between $\cH(p)$ and $\cH(p-1)$}

In this subsection we will construct recollements between $\cH(p)$ and $\cH(p-1)$ by using reduction functors and insertion functors.
First we give a useful result, which is of interest itself:

\begin{Lem}\label{adjointness for simples}
Let $\cC$ be an additive category, $S$ be an object in $\cC$ with $\End(S)$ a division ring. Denote by $\Delta:=\End(S)$ and $\widetilde{D}:=\Hom_{\Delta}(-, \Delta)$. Then there is an adjoint triple:
\[
\xymatrix@C=0.5cm{
\mod \Delta  \ar[rrr]^{-\otimes_{\Delta} S} &&& \cC.
\ar @/_1.5pc/[lll]_{\widetilde{D}\Hom_{\cC}(-, S)}
\ar @/^1.5pc/[lll]^{{\Hom_{\cC}(S, -)}}
 }
\]
\end{Lem}

\begin{pf}
The only non-trivial thing is
$(\widetilde{D}\Hom_{\cC}(-, S), -\otimes_{\Delta} S)$ form an adjoint functor. This follows from the natural isomorphisms
$$\aligned
 \Hom_{\Delta}(\widetilde{D}\Hom_{\cC}(C, S), \Delta^n) &\cong (\widetilde{D}(\widetilde{D}\Hom_{\cC}(C, S)))^n\\
 &\cong \Hom_{\cC}(C, S^n)  \\
 &\cong \Hom_{\cC}(C, \Delta^n \otimes_{\Delta} S),
 \endaligned
 $$
which are functorial in $C$ and $\Delta^n$ for any object $C\in\cC$ and $n\in\mathbb{N}$.
\end{pf}

The following lemma is a generalization of Lemma \ref{adjointness for simples}.

\begin{Lem}\label{adjointness for projectives}
Let $\cC$ be a Hom-finite abelian $\mathbf{k}$-linear category, let $\cal{S}$ be a full exact subcategory of $\cC$ with a projective generator $T$. Denote by $A:=\End_{\cC}(T)$. Then the following hold:
\begin{itemize}
  \item[(1)] the functor $\Hom_{\cC}(T,-)$ induces an equivalence
  $\xymatrix{\cS\ar^-{\sim}[r]& \mod A}$ with quasi-inverse $-\otimes_{A}T$;
  \item[(2)] the functor $-\otimes_{A}T: \mod A\to \cC$ admits a right adjoint $\Hom_{\cC}(T, -)$ and a left adjoint $D\Hom_{\cC}(-, DA\otimes_AT)$.
\end{itemize}
\end{Lem}

\begin{pf} The statement (1) follows from the assumption that $T$ is a projective generator. For statement (2), we
denote by
$q=D\Hom_{\cC}(-, DA\otimes_AT)$ and show that $(q, -\otimes_{A} T)$ forms an adjoint pair. The adjointness of $(-\otimes_{A} T, \Hom_{\cC}(T, -))$ follows from dual arguments.

Observe that for any $Y\in\cC$, there is a natural isomorphism:
\begin{equation}\label{natural equiv} \Hom_{A}(qY, DA)\cong D(qY) \cong  \Hom_{\cC}(Y, DA\otimes_A T).
\end{equation}
For any $M\in \mod A$, consider the following injective co-presentation of $M$:
$$\xymatrix{ 0\ar[r] & M\ar[r]& (DA)^n \ar^{f}[r]&(DA)^m}.$$
Now $T$ is a projective generator implies $-\otimes_{A}T$ is exact, yielding an exact sequence
$$\xymatrix{ 0\ar[r] & M\otimes_{A}T\ar[r]& (DA)^n\otimes_{A}T \ar^{f\otimes \id_{T}}[r]&(DA)^m\otimes_{A}T}.$$
This combining with (\ref{natural equiv}) yields the following commutative diagram:
$$\xymatrix{ 0\ar[r] & \Hom(qY, M) \ar[r]\ar^{\exists\,\eta_1}@{-->}[d]& \Hom(qY, (DA)^n) \ar^{\Hom(qY,f)}[r]\ar^{\eta_2}[d]\ar@{}|{\circlearrowleft}[rd]& \Hom(qY, (DA)^n) \ar^{\eta_3}[d]\\
0\ar[r] & \Hom(Y, M\otimes_{A}T)\ar[r]& \Hom(Y, (DA)^n\otimes_{A}T) \ar^{\Hom(Y,f\otimes \id_{T})}[r]&\Hom(Y,(DA)^m\otimes_{A}T),
}$$
where $\eta_2$ and $\eta_3$ are both natural isomorphisms. It follows that  $\eta_1$ is also a natural isomorphism, which claims the ajointness.
\end{pf}

Recall that $\cS_x$ is the set of simple sheaves in $\cal{H}$ concentrated in $x$, and $\tilde{\mathbf{S}}_{i}=\bigoplus_{S\in\cS_x}S_i$ for $i\in\bbZ/p\bbZ$.
%, where the indices are taken modulo $p$.
For any simple sheaf $S\in\cS_x$, $\End_{\cH(p)}(S_i)\cong\End_{\cH}(S):=\Delta$ is a finite skew field extension of $\mathbf{k}$. Then
$\tilde{\Delta}:=\End_{\cH(p)}(\tilde{\mathbf{S}}_{i})$ is a product of $n$-copies of $\Delta$, where $n$ is the order of the set $\cS_x$.

\begin{Prop}\label{recollement for H(p-1)} For any $i\in\bbZ$,
the following diagram is a recollement:
\begin{equation}\label{recollement graph of H(p-1)}
\xymatrix@C=0.5cm{
\mod\tilde{\Delta} \ar[rrr]^{-\otimes_{\tilde{\Delta}}\tilde{\mathbf{S}}_{p-i}} &&& \cH(p) \ar[rrr]^{\psi^i} \ar @/_1.5pc/[lll]_{{D}\Hom_{\cH(p)}(-, D\tilde{\Delta}\otimes_{\tilde{\Delta}}\tilde{\mathbf{S}}_{p-i})}  \ar
 @/^1.5pc/[lll]^{{\Hom_{\cH(p)}(\tilde{\mathbf{S}}_{p-i}, -)}} &&& \cH(p-1).
\ar @/_1.5pc/[lll]_{\psi_{i-1}} \ar
 @/^1.5pc/[lll]^{\psi_i}
 }
\end{equation}
\end{Prop}

\begin{pf}
Combining with Proposition \ref{adjoint sequence in abelian category}, Lemma \ref{exact of r and l} and Lemma \ref{ker and image}, there is an adjoint triple $\psi_i\dashv\psi^i\dashv \psi_{i+1}$, where $\psi_i$ and $\psi_{i+1}$ are fully faithful and $\Ker(\psi^i)=\add(\tilde{\mathbf{S}}_{p-i})=\im(-\otimes_{\tilde{\Delta}}\tilde{\mathbf{S}}_{p-i})$. Then by Lemma \ref{adjointness for projectives}, we obtain the required recollement.
\end{pf}

\subsection{Ladders between $D^b(\cH(p))$ and $D^b(\cH(p-1))$}

Given any abelian category $\cal{A}$, we denote by $D^b(\cal{A})$ its bounded derived category. An exact functor $F: \cal{A}\to \cal{B}$ between abelian categories extends to an exact functor
$F^{\ast}: D^b(\cal{A}) \to  $  $D^b(\cal{B})$. Note that the kernel of $F^{\ast}$ coincides with the full subcategory $D^b_{\Ker (F)}(\cal{A})$ consisting of the complexes in $D^b(\cal{A})$ with cohomology in $\Ker (F)$. In particular, if $\cal{A}$ is hereditary, then $\Ker (F^{\ast})=D^b(\Ker (F))$, and similarly $\im (F^{\ast})=D^b(\im (F))$.

The following result is helpful for lifting adjoint functors from abelian categories to their bounded derived categories.

\begin{Lem}\label{lifting from abelian to derive cat}\cite[Lem. 3.3.1]{[CK]} Let $F: \cal{A}\to\cal{B}$ be a fully faithful exact functor between abelian categories and suppose that $F$ admits an exact right adjoint $G: \cal{B}\to\cal{A}$.
\begin{itemize}
\item[(1)] $F^{\ast}$ and $G^{\ast}$ form an adjoint pair of exact functors and  $F^{\ast}$ is fully faithful.
\item[(2)] The inclusion $\Ker (G^{\ast}) \to D^b(\cal{B})$ admits a left adjoint which induces an equivalence  $D^b(\cal{B})/\im$$(G^{\ast})\cong\Ker (G^{\ast})$.
\end{itemize}
\end{Lem}

Observe that all the functors $\psi^i, \psi_i$ and $-\otimes_{\tilde{\Delta}}\tilde{\mathbf{S}}_{p-i}$ are exact for any $i\in \bbZ$. Hence they extend to exact functors between the bounded derive categories, which we denote by $\Psi^i, \Psi_i$ and  $T_i$ respectively. Then by Lemma \ref{lifting from abelian to derive cat}, $\Psi^i, \Psi_i$ are exact and $\Psi_i, T_i$ are fully faithful. These functors fit into an infinite ladder as follows.

\begin{Prop}\label{p ladder for derive category} There exists an infinite ladder of period $p$ as follows:
 \begin{equation}\label{infinite ladder}
\xymatrix@C=0.5cm{
D^b(\mod\tilde{\Delta})
\ar @/^-3pc/[rrr]_{\vdots}|{T_1} \ar @/^3pc/[rrr]^{\vdots}|{T_{-1}}
\ar[rrr]|{T_0}
&&& D^b(\cH(p))
\ar@/_1.5pc/[lll]|{T^{-1}}
\ar@/^1.5pc/[lll]|{T^{0}}
\ar[rrr]|{\Psi^0 }
\ar @/^-3pc/[rrr]_{\vdots}|{\Psi^1}
\ar @/^3pc/[rrr]^{\vdots}|{\Psi^{-1}}
&&& D^b(\cH(p-1)).
\ar@/_1.5pc/[lll]|{\Psi_{-1}}
\ar@/^1.5pc/[lll]|{\Psi_{0}}
}
\end{equation}
\end{Prop}

\begin{pf}
Firstly, according to Lemma \ref{lifting from abelian to derive cat}, the adjoint sequence (\ref{adjoint sequence of psi}) between $\cH(p)$ and $\cH(p-1)$ induces an infinite sequence of adjoint functors between their bounded derived categories:
\begin{equation}\label{adjoint sequence of Psi} \cdots\dashv\Psi^{-1}\dashv\Psi_{-1}\dashv\Psi^0\dashv\Psi_0\dashv\Psi^1\dashv \cdots.
\end{equation}
Moreover, for each $i\in\bbZ$, $\psi^i\psi_i=\id_{\cH(p-1)}$ implies $\Psi^i\Psi_i=\id_{D^b(\cH(p-1))}$. It follows that $\Psi^{i}$ is a quotient functor and $$\Ker(\Psi^i)=D^b(\Ker(\psi^i))=D^b(\im(-\otimes_{\tilde{\Delta}}\tilde{\mathbf{S}}_{p-i}))=\im(T_i).$$
Consequently, there is an exact sequence of triangulated categories:
\[
\xymatrix@C=0.5cm{
D^b(\mod\tilde{\Delta})  \ar[rrr]^{T_i}
&&& D^b(\cH(p))
\ar[rrr]^{\Psi^i }
&&& D^b(\cH(p-1)).
}
 \]
This induces a recollement by \cite[Thm. 2.1]{[CPS2]}:
\[
\xymatrix@C=0.5cm{
D^b(\mod\tilde{\Delta})  \ar[rrr]|{T_i}
&&& D^b(\cH(p))
\ar@/_1.5pc/[lll]|{}
\ar@/^1.5pc/[lll]|{}
\ar[rrr]|{\Psi^i }
&&& D^b(\cH(p-1)).
\ar@/_1.5pc/[lll]|{\Psi_{i-1}}
\ar@/^1.5pc/[lll]|{\Psi_{i}}
}
 \]
Denote the left and right adjoint functors of $T_i$ by $T_{l}^i$ and $T_{r}^i$ respectively. We claim that $T_{r}^i=T_{l}^{i+1}$ for each $i\in\bbZ$.
In fact, by \cite[Thm. 1.1]{[CPS1]} the left recollement
\[
\xymatrix@C=0.5cm{
D^b(\cH(p-1)) \ar[rrr]|{\Psi_{i}}
&&& D^b(\cH(p))
\ar@/_1.5pc/[lll]|{\Psi^{i}}
\ar@/^1.5pc/[lll]|{\Psi^{i+1}}
}
 \] can be completed to a recollement, which has the following form by Lemma \ref{lifting from abelian to derive cat} and its dual:
 \[
\xymatrix@C=0.5cm{
D^b(\cH(p-1)) \ar[rrr]|{\Psi_{i}}
&&& D^b(\cH(p))
\ar@/_1.5pc/[lll]|{\Psi^{i}}
\ar@/^1.5pc/[lll]|{\Psi^{i+1}}
\ar[rrr]|{ }
&&& D^b(\mod\tilde{\Delta}).
\ar@/_1.5pc/[lll]|{T_{i}}
\ar@/^1.5pc/[lll]|{T_{i+1}}
}
 \]
  Here we use the following equivalences
 $$D^b(\cH(p))/\im(\Psi_i)\cong\Ker(\Psi^{i+1})=\im(T_{i+1})\cong D^b(\mod\tilde{\Delta}).$$
 It follows that the right adjoint functor of $T_i$ coincides with the left adjoint functor of $T_{i+1}$, i.e, $T_{r}^i=T_{l}^{i+1}$, which we denote by $T^i$. Then the infinite ladder (\ref{infinite ladder}) follows.

Now we consider the periodicity of the ladder (\ref{infinite ladder}). Recall that there is a one-to-one correspondence between equivalence classes of recollements of $D^b(\cH(p))$ and TTF-triples of $D^b(\cH(p))$, see for example \cite[Prop. 2.1]{[AKLY]}. Observe that $\Ker(\Psi^{np+i})=\Ker(\Psi^{i})$ and
$\im(\Psi_{np+i})=\im(\Psi_{i})$ for any $n\in\bbZ$ and $0\leq i\leq p-1$. Hence there is an unbounded TTF-tuple which is periodic of period $p$:
$$(\cdots, \im(\Psi_{i-1}), \Ker(\Psi^{i}), \im(\Psi_{i}), \Ker(\Psi^{i+1}), \im(\Psi_{i+1}),\cdots ).$$
It follows that the ladder (\ref{infinite ladder}) is periodic of period $p$.
\end{pf}

\subsection{Composition of recollements}
In this subsection, we fix an integer $q$ with $0<q<p$ to describe the recollements between the categories $\cH(p)$ with $\cH(p-q)$.

We use the notation $\underline{q}$ to denote an arbitrary sequence $(i_1, i_2,\cdots,i_q)$ of integers with
$0\leq i_1<i_2<\cdots<i_{q}<p$. Set $\underline{1}=(1,1,\cdots 1)$ with $q$-entries, and define $\underline{q}\pm \underline{1}=(i_1\pm 1, i_2\pm 1 ,\cdots,i_q\pm 1)$ via componentwise addition.
We define the functor
$\psi^{\underline{q}}:=\psi^{i_1}\psi^{i_2}\cdots \psi^{i_{q-1}}\psi^{i_q}$ as the
composition of the following functors:
\[
\xymatrix@C=0.5cm{
\cH(p) \ar[rr]^-{\psi^{i_q}}
 &&
 \cH(p-1) \ar[rr]^-{\psi^{i_{q-1}}}&&
 \cdots\ar[rr]^-{\psi^{i_2}}
 &&
 \cH(p-q+1) \ar[rr]^-{\psi^{i_1}}
 && \cH(p-q).  }
\]
Similarly, we define
$\psi_{\underline{q}}=\psi_{i_q}\psi_{i_{q-1}}\cdots \psi_{i_2}\psi_{i_1}:\cH(p-q)\to\cH(p).$
As a consequence of Proposition \ref{adjoint sequence in abelian category}, we have an infinite sequence of adjoint functors
\begin{equation}\label{infinite seq of comp}\cdots\dashv\psi_{\underline{q}-\underline{1}}\dashv \psi^{\underline{q}}\dashv\psi_{\underline{q}}\dashv \psi^{\underline{q}+\underline{1}}\dashv \psi_{\underline{q}+\underline{1}}\dashv \cdots.
\end{equation}

Let $\cS(\underline{q})$ be the full subcategory of $\cH(p)$ generated by the simple objects $S_{i_j}$'s, where $S\in\cS_x$ and $1\leq j\leq q$. Then $\cS(\underline{q})$ is equivalent to a module category $\mod A(\underline{q})$, where $A(\underline{q})$ is a product of path algebras of Dynkin type $\mathbb{A}$$_n$'s (determined by $\underline{q}$). Denote by ${P_{\underline{q}}}$ the projective generator of $\cS(\underline{q})$.

Now we can obtain a recollement between $\cH(p)$ and $\cH(p-q)$ as follows, which can be thought as certain composition of recollements of the form (\ref{recollement graph of H(p-1)}).
 \begin{Prop}\label{recollement for general q}
For any $0< q< p$, the following diagram is a recollement:
 \[
\xymatrix@C=0.5cm{
\mod A(\underline{q}) \ar[rrr]^{-\otimes {P_{\underline{q}}}} &&& \cH(p) \ar[rrr]^{\psi^{\underline{q}}} \ar @/_1.5pc/[lll]_{D\Hom(-, DA(\underline{q})\otimes{P_{\underline{q}}})}  \ar
 @/^1.5pc/[lll]^{{\Hom({P_{\underline{q}}}, -)}} &&& \cH(p-q).
\ar @/_1.5pc/[lll]_{\psi_{\underline{q}-\underline{1}}} \ar
 @/^1.5pc/[lll]^{\psi_{\underline{q}}}
 }
\]

\end{Prop}

\begin{pf} 
By (\ref{infinite seq of comp}) there is an adjoint triple $\psi_{\underline{q}-\underline{1}}\dashv\psi^{\underline{q}}\dashv \psi_{\underline{q}}$. Moreover, $\psi_{\underline{q}-\underline{1}}$ (resp. $\psi_{\underline{q}-\underline{1}}$) is a composition of fully faithful functors, hence itself is fully faithful. On the other hand, Lemma \ref{adjointness for projectives} yields the adjointness of the left half recollement and that $-\otimes {P_{\underline{q}}}$ is fully faithful. Hence it remains to show that $\Ker(\psi^{\underline{q}})=\im(-\otimes {P_{\underline{q}}})$.
For this let $E$ be an indecomposable object in $\cH(p)$ of the form (\ref{p-cycle}). By definition, $\psi^{\underline{q}}(E)=0$ if and only if $E_{j}=0$ for any $j\notin\{i_1, i_2,\cdots, i_q\}$, i.e, $E$ belongs to the subcategory $\cS(\underline{q}).$ It follows that $\Ker(\psi^{\underline{q}})=\cS(\underline{q})=\im(-\otimes {P_{\underline{q}}})$. This finishes the proof.
\end{pf}

Passing to the bounded derived categories, we denote by
$\Psi^{\underline{q}}=\Psi^{i_1}\Psi^{i_2}\cdots \Psi^{i_q}$ and
$\Psi_{\underline{q}}=\Psi_{i_q}\cdots \Psi_{i_2}\Psi_{i_1}$ in a similar way as before.
Let $T_{\underline{q}}$ be the derived functor of $-\otimes {P_{\underline{q}}}$. By similar arguments as in the proof of Proposition \ref{p ladder for derive category}, we obtain that $T_{\underline{q}}$ admits a right adjoint $T^{\underline{q}}$, and these functors fit into an infinite ladder:

 \begin{Prop}\label{ladders for general H} %For each $0\leq i_1<i_2<\cdots< i_q<p$,
 The following diagram is an infinite ladder of period $p$:
\[
\xymatrix@C=0.5cm{
D^b(\mod A(\underline{q}))
\ar@/^-3pc/[rrr]_{\vdots}|{T{\underline{q}+\underline{1}}} \ar @/^3pc/[rrr]^{\vdots}|{T_{\underline{q}-\underline{1}}}
\ar[rrr]|{T_{\underline{q}}}
&&& D^b(\cH(p))
\ar@/_1.5pc/[lll]|{T^{\underline{q}-\underline{1}}}
\ar@/^1.5pc/[lll]|{T^{\underline{q}}}
\ar[rrr]|{\Psi^{\underline{q}}}
\ar@/^-3pc/[rrr]_{\vdots}|{\Psi^{\underline{q}+\underline{1}}}
\ar @/^3pc/[rrr]^{\vdots}|{\Psi^{\underline{q}-\underline{1}}}
&&& D^b(\cH(p-q)).
\ar@/_1.5pc/[lll]|{\Psi_{\underline{q}-\underline{1}}}
\ar@/^1.5pc/[lll]|{\Psi_{\underline{q}}}
}
\]
\end{Prop}

\section{Recollements and Ladders for weighted projective lines}

In this section we focus on the study of recollements for weighted projective lines in three different levels, including the category of coherent sheaves, its bounded derived category and also the stable category of vector bundles.

\subsection{$p$-cycle interpretation for weighted projective lines}
We first recall the alternative description of weighted projective lines in the sense of Geigle-Lenzing \cite{[GL1]} by $p$-cycle construction due to Lenzing \cite{[L]}.

Let $\bbX_0$ be the usual projective line over $\mathbf{k}$, and let $\blambda=(\lambda_1,\lambda_2,\cdots,\lambda_t)$ be a sequence of pairwise distinct points from $\bbX_0$ and ${\mathbf{p}}=(p_1, p_2,\cdots, p_t)$ be a sequence of positive integers. Associated to these data there is a weighted projective line $\bbX({\mathbf{p}}, \blambda)$ with the weight sequence ${\mathbf{p}}$ and parameter sequence $\blambda$, we refer to \cite{[GL1]} for details.

Let inductively $\bbX_i$ be the exceptional curve obtained from $\bbX_{i-1}$ by inserting weight $p_i$ in $\lambda_{i}$, i.e, forming the category of $p_i$-cycles in $\coh\bbX_{i-1}$ which are concentrated in $\lambda_i$. Then $\bbX_t$ is isomorphic to the weighted projective line $\bbX({\mathbf{p}}, \blambda)$ (this statement has been shown in \cite[Section 4]{[L]}, an alternative proof is given by Hubery \cite{Hubery} recently). In the following we identify $\bbX_t$ with $\bbX({\mathbf{p}}, \blambda)$.

From this construction, we see that for any $1\leq i, j\leq t$, if $j\leq i-1$, then there is a unique simple sheaf in $\bbX_j$ concentrated in $\lambda_i$, which will be denoted by $S_{\lambda_i}$; if $j\geq i$, then there are $p_i$-simples in $\bbX_j$ concentrated in $\lambda_i$, which will be denoted by $S_{\lambda_i,j}$'s, where the index $j\in\bbZ$ is taken modulo $p_i$.

\subsection{Recollements of $\coh\bbX({\mathbf{p}}, \blambda)$}

In this subsection, we aim to classify the recollements of the category $\coh\bbX({\mathbf{p}}, \blambda)$.

Let ${\mathbf{q}}=(q_1, q_2,\cdots, q_t)$ be a sequence of integers with $0\leq q_i\leq p_i-1$ for each $i$, which will be simplified denoted by ${\mathbf{q}}\prec{\mathbf{p}}$ below. For $0\leq i\leq t$, we denote by $$\bbY_i=\bbX(p_1,\cdots, p_{i-1}, p_i, p_{i+1}-q_{i+1}, \cdots, p_{t}-q_{t};\; \lambda_1,\cdots, \lambda_t)$$ and $$\hat{\bbY}_i=\bbX(p_1,\cdots, p_{i-1}, 1, p_{i+1}-q_{i+1}, \cdots, p_{t}-q_{t};\; \lambda_1,\cdots, \lambda_t).$$ Then we have $\bbY_0=\hat{\bbY}_0=\bbX({\mathbf{p}}-{\mathbf{q}}, \blambda)$, $\bbY_t=\bbX({\mathbf{p}}, \blambda)$ and $\bbY_i$ (resp. $\bbY_{i-1}$) is obtained from $\hat{\bbY}_i$ by $p_i$-cycle (resp. $(p_{i}-q_{i})$-cycle) construction for $1\leq i\leq t$. By Proposition \ref{recollement for general q}, there is an adjoint triple
 \[
\xymatrix@C=0.5cm{
\coh{\bbY}_i=(\coh\hat{\bbY}_i)(p_i;  \lambda) \ar[rrr]^{\psi^{\underline{q_i}}_{{}}}
 &&& (\coh\hat{\bbY}_i)(p_i-q_i; \lambda)=\coh{\bbY}_{i-1}
\ar @/_1.5pc/[lll]_{\psi^{ }_{\underline{q_i}-{\underline{1}}}} \ar
 @/^1.5pc/[lll]^{\psi^{ }_{\underline{q_i}}}
 }
\]
for any
integer sequence $\underline{q_i}=(i_1, i_2,\cdots, i_{q_i})$ with $0\leq i_1<i_2<\cdots< i_{q_i}<p_i$.
In this case, $\Ker(\psi_{ }^{\underline{q_i}})$ is a Serre subcategory of $\coh{\bbY}_i$ generated by simple sheaves $S_{\lambda_i,i_j}\ (\; 1\leq j\leq q_i)$, which is equivalent to the module category $\mod A(\underline{q_i})$.

For any sequence $\underline{{\mathbf{q}}}=(\underline{q_1}, \underline{q_2},\cdots, \underline{q_t})$, we denote by $\psi^{\underline{{\mathbf{q}}}}:=
\psi^{\underline{q_1}}_{{}}\cdots\psi^{\underline{q_{t-1}}}_{{}}\psi^{\underline{q_t}}_{{}}$ the following composition of functors \[
\xymatrix@C=0.5cm{
\coh\bbX({\mathbf{p}}, \blambda)=\coh{\bbY}_t \ar[rr]^-{\psi^{\underline{q_t}}_{{}}}
 &&
 \coh{\bbY}_{t-1} \ar[rr]^-{\psi^{\underline{q_{t-1}} }_{{}}}
 &&
 \cdots \ar[rr]^-{\psi^{\underline{q_{1}} }_{{}}}
 &&
 \coh{\bbY}_{0} =\coh\bbX({\mathbf{p}}-{\mathbf{q}}, \blambda),
 }
\]
and similarly denote by $$\psi_{\underline{{\mathbf{q}}}}:=\psi_{\underline{q_t}}^{{}}
\psi_{\underline{q_{t-1}}}^{{}}\cdots\psi_{\underline{q_1}}^{{}}: \ \coh\bbX({\mathbf{p}}-{\mathbf{q}}, \blambda)\to \coh\bbX({\mathbf{p}}, \blambda).$$
Observe that $\Ker(\psi^{\underline{{\mathbf{q}}}})$ is a Serre subcategory of $\coh{\bbY}_t$ generated by simple sheaves $S_{\lambda_i,i_j}$'s with $1\leq i\leq t$ and $1\leq j\leq q_i$, which we denote by $\cS(\underline{{\mathbf{q}}})$. Then $\cS(\underline{{\mathbf{q}}})$ is equivalent to the module category $\mod A(\underline{{\mathbf{q}}})$, where $A(\underline{{\mathbf{q}}})=A(\underline{q_1})\times\cdots \times A(\underline{q_t})$.
Denote by $P_{\underline{{\mathbf{q}}}}$ the projective generator of $\cS(\underline{{\mathbf{q}}})$. We have the following main result:

\begin{Thm}\label{theorem of recollement for coh} Up to equivalence, each recollement of $\coh\bbX({\mathbf{p}}, \blambda)$ has the following form for some integer sequence ${\mathbf{q}}\prec{\mathbf{p}}$ and some sequence $\underline{{\mathbf{q}}}$:

\begin{equation}\label{recollement for coh}
\xymatrix@C=0.5cm{
\mod A(\underline{{\mathbf{q}}}) \ar[rrr]^{-\otimes{P_{\underline{{\mathbf{q}}}}}}
&&& \coh\bbX({\mathbf{p}}, \blambda) \ar[rrr]^{\psi^{\underline{{\mathbf{q}}}}} \ar @/_1.5pc/[lll]_{D\Hom(-, DA(\underline{{\mathbf{q}}})\otimes{P_{\underline{{\mathbf{q}}}}})}  \ar
 @/^1.5pc/[lll]^{{\Hom({P_{\underline{{\mathbf{q}}}}}, -)}}
 &&& \coh\bbX({\mathbf{p}}-{\mathbf{q}}, \blambda).
\ar @/_1.5pc/[lll]_{\psi_{\underline{{\mathbf{q}}}-\underline{{\mathbf{1}}}}} \ar
 @/^1.5pc/[lll]^{\psi_{\underline{{\mathbf{q}}}}}
 }
\end{equation}
 \end{Thm}

\begin{pf}
Using similar proof as in Proposition \ref{recollement for general q}, one shows (\ref{recollement for coh}) is a recollement. Recall that there is a bijection between equivalent classes of recollements of $\coh\bbX$ and bilocalizing Serre subcategories of $\coh\bbX$, see for example \cite[Thm. 4.3]{[PV]}. Observe that the corresponding bilocalizing Serre subcategories of (\ref{recollement for coh}) is just $\cS(\underline{{\mathbf{q}}})$. Hence it suffices to show that any bilocalizing Serre subcategories of $\coh\bbX$ has the form $\cS(\underline{{\mathbf{q}}})$ for some $\underline{{\mathbf{q}}}$, this follows from \cite[Prop. 9.2]{[GL]} and its dual.
\end{pf}

\subsection{Ladders of $D^b(\coh\bbX({\mathbf{p}}, \blambda))$}

Now we consider the ladders of $D^b(\coh\bbX({\mathbf{p}}, \blambda))$. By Lemma \ref{lifting from abelian to derive cat} we obtain a family of derived functors $\Psi^{\underline{{\mathbf{q}}}}, \Psi_{\underline{{\mathbf{q}}}}$ and $T_{\underline{{\mathbf{q}}}}$ from the exact functors $\psi^{\underline{{\mathbf{q}}}}, \psi_{\underline{{\mathbf{q}}}}$ and $-\otimes{P_{\underline{{\mathbf{q}}}}}$ respectively. Similarly as in Proposition \ref{ladders for general H}, we have
$\Psi^{\underline{{\mathbf{q}}}}=
\Psi^{\underline{q_1}}_{{}}\cdots\Psi^{\underline{q_{t-1}}}_{{}}
\Psi^{\underline{q_t}}_{{}}$ and $\Psi_{\underline{{\mathbf{q}}}}=\Psi_{\underline{q_t}}^{{}}
\Psi_{\underline{q_{t-1}}}^{{}}\cdots\Psi_{\underline{q_1}}^{{}}$, and $T_{\underline{{\mathbf{q}}}}$ admits a right adjoint functor $T^{\underline{{\mathbf{q}}}}$, and all these functors fit into an infinite ladder: 

 \begin{Thm}\label{main theorem for der cat} 
 For any sequences ${\mathbf{q}}$ and $\underline{{\mathbf{q}}}$ as in Theorem
 \ref{theorem of recollement for coh},
  there is an infinite ladder as follows, which is periodic of period l.c.m.$(p_1,p_2,\cdots,p_t)$:
\[
\xymatrix@C=0.5cm{
D^b(\mod A(\underline{{\mathbf{q}})})
\ar@/^-3pc/[rrr]_{\vdots}|{T_{\underline{{\mathbf{q}}}+\underline{{\mathbf{1}}}}}
\ar @/^3pc/[rrr]^{\vdots}|{T_{\underline{{\mathbf{q}}}-\underline{{\mathbf{1}}}}}
\ar[rrr]|{T_{\underline{{\mathbf{q}}}}}
&&& D^b(\coh\bbX({\mathbf{p}}, \blambda))
\ar@/_1.5pc/[lll]|{T^{\underline{{\mathbf{q}}}-\underline{{\mathbf{1}}}}}
\ar@/^1.5pc/[lll]|{T^{\underline{{\mathbf{q}}}}}
\ar[rrr]|{\Psi^{\underline{{\mathbf{q}}}}}
\ar@/^-3pc/[rrr]_{\vdots}|{\Psi^{\underline{{\mathbf{q}}}+\underline{{\mathbf{1}}}}}
\ar @/^3pc/[rrr]^{\vdots}|{\Psi^{\underline{{\mathbf{q}}}-\underline{{\mathbf{1}}}}}
&&& D^b(\coh\bbX({\mathbf{p}}-{\mathbf{q}}, \blambda)).
\ar@/_1.5pc/[lll]|{\Psi_{\underline{{\mathbf{q}}}-\underline{{\mathbf{1}}}}}
\ar@/^1.5pc/[lll]|{\Psi_{\underline{{\mathbf{q}}}}}
}
 \]
\end{Thm}

\subsection{Ladders for the stable categories of vector bundles}

Let $\bbX_t=\bbX({\mathbf{p}}, \blambda)$ be a weighted projective line of weight type ${\mathbf{p}}=(p_1, p_2,\cdots,p_t)$ with parameter data $\blambda=(\lambda_1,\lambda_2\cdots,\lambda_t)$.
Denote by $\vect\bbX_t$ the full subcategory of $\coh\bbX_t$ formed by all torsion-free sheaves (i.e. vector bundles). We introduce briefly a Frobenius exact structure on $\vect\bbX_t$, for more details we refer to \cite[Sect. 3]{[KLM]}.

A sequence $\xi: 0\to X'\to X\to X''\to 0$ in $\vect\bbX_t$ is called \emph{distinguished exact} if for each line bundle $L$ the induced sequence $$\Hom(L, \xi): \quad 0\to \Hom(L, X')\to \Hom(L, X)\to \Hom(L, X'')\to 0$$ is exact. By using Serre duality, we know that $\xi$ is distinguished exact if and only if $\Hom(\xi, L)$ is exact for each line bundle $L$.
According to \cite[Prop. 3.2]{[KLM]}, the distinguished exact sequences define a Frobenius exact structure on $\vect\bbX_t$, such that the indecomposable projectives (resp. injectives) are exactly the line bundles. It follows that the associated stable category $\underline{\vect}\bbX_t$ obtained from $\vect\bbX_t$ by factoring out all the line bundles is a triangulated category.

For simplicity, we denote by
$\bbY:=\bbX_{t-1}=\bbX(p_1, p_2,\cdots, p_{t-1}; \lambda_1,\lambda_2, \cdots,\lambda_{t-1})$, and use the notation $\bbY(q):=\bbX(p_1, p_2,\cdots, p_{t-1},
q; \lambda_1, \lambda_2,\cdots,\lambda_t)$ for any $1\leq q\leq p:=p_t$.
It follows that $\coh \bbY(q)=(\coh\bbY)(q; \lambda_{t})$. Then by Proposition \ref{adjoint sequence in abelian category}, there exists an infinite sequence of adjoint functors between $\coh\bbY(p)$ and $\coh\bbY(p-1)$:
$$ \cdots\dashv\psi^{-1}\dashv\psi_{-1}\dashv\psi^0\dashv\psi_0\dashv\psi^1\dashv \cdots  .
$$
Observe that for any $i\in\bbZ$, the functor
$\psi^i: \coh\bbY(p)\to \coh\bbY(p-1)$
preserves the rank of coherent sheaves, in particular, it preserves torsion sheaves (rank$=0$) and vector bundles (rank$>0$) respectively, hence it restricts to a functor between the subcategories of vector bundles, which we still denote by $\psi^i: \vect\bbY(p)\to \vect\bbY(p-1)$. Similarly, we have a restriction functor $\psi_i: \vect\bbY(p-1)\to \vect\bbY(p)$.

\begin{Lem}
For any $i\in\bbZ$, the restriction functors $\psi^i$ and $\psi_i$ on vector bundles preserve distinguished exact sequences.
\end{Lem}

\begin{pf}We use the adjointness $\psi^i\dashv \psi_i$ to show that $\psi_i:\vect\bbY(p-1)\to \vect\bbY(p)$ preserves distinguished exact sequences.
The proof for $\psi^i$ is similar.

Let $\xi: 0\to X'\to X\to X''\to 0$ be a distinguished exact sequence in $\vect\bbY(p-1)$. For any line bundle $L\in\vect\bbY(p)$, we know that $\psi^i(L)$ is a line bundle in $\vect\bbY(p-1)$. Hence by definition we get an exact sequence
$$0\to \Hom(\psi^i(L), X')\to \Hom(\psi^i(L), X)\to \Hom(\psi^i(L), X'')\to 0.$$
Now using the adjointness $\psi^i\dashv\psi_i$, we obtain an exact sequence
$$0\to \Hom(L, \psi_i(X'))\to \Hom(L, \psi_i(X))\to \Hom(L, \psi_i(X''))\to 0, $$ which means $\psi_i(\xi)$ is distinguished exact. We are done.
\end{pf}

Consequently, for any integer $i$, the functors $\psi^i$ and $\psi_i$ induce triangle functors between the stable categories:
$$\bar{\psi}^i: \underline{\vect}\bbY(p) \to\underline{\vect}\bbY(p-1)\quad\text{and}\quad
\bar{\psi}_i: \underline{\vect}\bbY(p-1) \to\underline{\vect}\bbY(p).$$

\begin{Prop}\label{adjoint tuples for stable category} For any $i\in\bbZ$, we have $\bar{\psi}^i\bar{\psi}_i=\id_{\underline{\vect}\bbY(p-1))}=\bar{\psi}^i\bar{\psi}_{i-1}$ and $(\bar{\psi}_{i-1}, \bar{\psi}^i, \bar{\psi}_i)$ form an adjoint triple. Consequently, $\bar{\psi}_i$ is fully faithful and $\bar{\psi}^i$ is a quotient functor, and there exists an infinite sequence of adjoint functors
\begin{equation}\label{adjoint sequence in stab}\cdots\dashv \bar{\psi}^{-1}\dashv\bar{\psi}_{-1} \dashv \bar{\psi}^0 \dashv \bar{\psi}_0\dashv \bar{\psi}^{1} \dashv\cdots.
\end{equation}
%For each $0\leq i\leq p-1$, there exists a sequence of adjoint pairs between the stable categories of vector bundles $\underline{\vect}\bbY(p)$ and $\underline{\vect}\bbY(p-1)$:
%$$(\bar{\psi}^0, \bar{\psi}_0, \cdots, \bar{\psi}^i, \bar{\psi}_i, \cdots, \bar{\psi}^{p-1}, \bar{\psi}_{p-1}, \bar{\psi}^0).$$
\end{Prop}

\begin{pf} These follow from Proposition \ref{adjoint sequence in abelian category} and Lemma \ref{exact of r and l}.
 \end{pf}

The following lemma plays a key role for our further discussion.
\begin{Lem}\label{kernel for bar of ri}
For any $i\in\bbZ$, $\Ker(\bar{\psi^i})$ is the subcategory of $\underline{\vect}\bbY(p)$ consisting of $p$-cycles of the form
 $$\xymatrix{
E_0\ar@{=}[r] &  \cdots \ar@{=}[r] & E_{i-1}\ar[r]^{x_{i-1}} & E_{i}\ar[r]^{x_{i}} & E_{i+1}\ar@{=}[r] &  \cdots \ar@{=}[r] & E_{0}(x),
}
$$ where
$E_{0}$ is a direct sum of line bundles in $\coh\bbY$.
\end{Lem}

\begin{pf}
Observe that $\Ker(\bar{\psi}^{np+i})=\Ker(\bar{\psi^i})$ for any $n\in\bbZ$, hence we can restrict to $0\leq i\leq p-1$. Fix such an integer $i$ we let $E\in \Ker(\bar{\psi^i})$ be an indecomposable vector bundle over $\bbY(p)$.  Assume $E$ has the following $p$-cycle form: %(\ref{p-cycle})
\begin{equation}\label{exp of E in the proof}\xymatrix@C=0.5cm{
E_0\ar[r]^{x_0} & \cdots \ar[r]& E_{i-1}\ar[r]^{x_{i-1}}& E_{i}\ar[r]^{x_{i}} & E_{i+1}\ar[r]^{x_{i+1}}& \cdots \ar[r]&  E_{p-1}\ar[r]^{x_{p-1}} & E_{0}(x),
}\end{equation}
where each $E_i\in\vect\bbY$. Then $\bar{\psi}^i(E)=0$ means $$\psi^i(E):=\xymatrix@C=0.5cm{
E_0\ar[r]^{x_0} & \cdots \ar[r]& E_{i-1}\ar[rr]^{x_ix_{i-1}}&& E_{i+1}\ar[r]^{x_{i+1}}& \cdots \ar[r]&  E_{p-1}\ar[r]^{x_{p-1}} & E_{0}(x)
}$$ decomposes as a direct sum of line bundles, say, $\psi^i(E)=\bigoplus_{j}L_{j}$, where $L_{j}$ is a line bundle in $\vect\bbY(p-1)$. Suppose $L_j$ has the following $(p-1)$-cycle form:
$$\xymatrix{
L_{j,0}\ar[r]^{x'_{j, 0}} & \cdots \ar[r]& L_{j, i-1}\ar[r]^{x'_{j, i-1}} & L_{j, i+1}\ar[r]^{x'_{j, i+1}}& \cdots \ar[r]&  L_{j, p-1}\ar[r]^{x'_{j, p-1}} & L_{j, 0}(x)
}.$$
Then for each $0\leq k\leq p-1$ with $k\neq i-1$ or $i$, we have $E_k=\bigoplus_j L_{j,k}$ and $x_k={\rm diag}(x_{j, k}')_j$.

By \cite[Prop.2.6]{[P1]}, the recollement (\ref{recollement graph of H(p-1)}) yields a canonical exact sequence in $\coh\bbY(p)$:
$$\xymatrix{
0 \ar[r] & \Ker(\mu_E) \ar[r] & \psi_{i-1}\psi^i(E) \ar[r]^-{\mu_E} & E \ar[r]&  \Coker(\mu_E) \ar[r]  & 0,
}
$$
where $\Ker(\mu_E)$ and $\Coker(\mu_E)$
belong to the subcategory %$\im(-\otimes_\Delta S_{p-i}):=\add{S_{p-i}}$
of $\coh\bbY(p)$ generated by $S_{p-i}$. Note that $\psi_{i-1}\psi^i(E)$ is a vector bundle and there are no non-zero homomorphisms from torsion sheaves to vector bundles. It follows that $\Ker(\mu_E)=0$, hence there is a short exact sequence as follows for some $n\in\mathbb{N}$:
$$\xymatrix{
0 \ar[r]  & \psi_{i-1}\psi^i(E) \ar[r]^-{\mu_E} & E \ar[r]&  S_{p-i}^{\oplus n} \ar[r]  & 0.
}
$$
Now $E$ is indecomposable implies $\Ext^1(S_{p-i}, \psi_{i-1}L_{j})\neq 0$ for each line bundle $L_{j}$. Then by Proposition \ref{extensions and homs}, $L_{j}$ must have the form
$$\xymatrix{
L_{j,0}\ar@{=}[r] & \cdots \ar@{=}[r]& L_{j, i-1}\ar[r]^{x'_{j, i-1}} & L_{j, i+1}\ar@{=}[r]& \cdots \ar@{=}[r]&  L_{j, p-1}\ar@{=}[r]& L_{j, 0}(x)
}.$$
It follows that each morphism $x_k$ in (\ref{exp of E in the proof}) is an isomorphism for $k\neq i-1$ or $i$. This finishes the proof.
\end{pf}

\begin{Rem} Let
$\cD$ be the full subcategory of $\underline{\vect}\bbY(2)$ consisting of the 2-cycles $E_0\to E_1\to E_0(x)$ with $E_0, E_1\in\vect\bbY$ and $E_0$ being a direct sum of line bundles.

(1) If we take $p=2$ and $i=1$ in the above lemma, we obtain that $\cD=\Ker(\psi^1)$, which is a triangulated subcategory of  $\underline{\vect}\bbY(2)$.

(2) If we take $t=3$, i.e, $\bbY$ has weight type $(p_1, p_2)$, then each indecomposable vector bundle in $\coh\bbY$ is a line bundle, hence $\cD=\underline{\vect}\bbY(2)$.
\end{Rem}

\begin{Prop} The following diagram is an infinite
 ladder of period $p$: \[
\xymatrix@C=0.5cm{
\cD  \ar[rrr]|{\bar{t}_{1}}%|{\bar{t}_0}
\ar@/_3pc/[rrr]_{\vdots}%|{\bar{t}_{1}}
\ar@/^3pc/[rrr]^{\vdots}%|{\bar{t}_{-1}}
&&& \underline{\vect}\bbY(p)
\ar@/_3pc/[rrr]_{\vdots}|{\bar{\psi}^{2}}
\ar@/^3pc/[rrr]^{\vdots}|{\bar{\psi}^{0}}
\ar@/_1.5pc/[lll]
\ar@/^1.5pc/[lll]
\ar[rrr]|{\bar{\psi}^1}
&&& \underline{\vect}\bbY(p-1),
\ar@/_1.5pc/[lll]|{\bar{\psi}_{0}}
\ar@/^1.5pc/[lll]|{\bar{\psi}_1}
}
 \]
where the functor
$\bar{t}_1=\bar{\psi}_{p-1}\bar{\psi}_{p-2}\cdots\bar{\psi}_{2}|_{\cD}$.
\end{Prop}

\begin{pf}
By Proposition \ref{adjoint tuples for stable category}, there is an adjoint triple $\bar{\psi}_{0}\dashv\bar{\psi}^{1}\dashv\bar{\psi}_{1}$, where $\bar{\psi}_{0}$ and $\bar{\psi}_{1}$ are fully faithful, and $\bar{\psi}^{1}$ is a quotient functor. By definition $\bar{t}_1:\cD\to\underline{\vect}\bbY(p)$ sends an object
$\xymatrix{E_0\ar[r]^{x_0}& E_1\ar[r]^{x_1}& E_0(x)}$ to $$\xymatrix{E_0\ar[r]^{x_0}& E_1\ar[r]^{x_1}& E_2\ar@{=}[r]& \cdots \ar@{=}[r]
& E_{p-1}\ar@{=}[r] & E_0(x)}.$$
Then by Lemma \ref{kernel for bar of ri} we see that 
$\Ker(\bar{\psi}^{1})=\im(\bar{t}_1)$. Hence we obtain a recollement by \cite[Thm. 2.1]{[CPS2]} as follows:
\begin{equation}\label{recollement diagram of D}
\xymatrix@C=0.5cm{
\cD  \ar[rrr]|{\bar{t}_{1}}%|{\bar{t}_0}
&&& \underline{\vect}\bbY(p)
\ar@/_1.5pc/[lll]
\ar@/^1.5pc/[lll]
\ar[rrr]|{\bar{\psi}^1}
&&& \underline{\vect}\bbY(p-1),
\ar@/_1.5pc/[lll]|{\bar{\psi}_{0}}
\ar@/^1.5pc/[lll]|{\bar{\psi}_1}
}
\end{equation}
Note that the category $\underline{\vect}\bbY(p)$ admits Serre duality, by using the reflecting approach introduced by J{\o}gensen in \cite{[J]}, we finally obtain a ladder, which should have the form (\ref{ladder diagram of stab cate}) by the exactness sequence (\ref{adjoint sequence in stab}).

Recall that there is a bijection between recollements of $\underline{\vect}\bbY(p)$ and triangle TTF-triples in $\underline{\vect}\bbY(p)$. In particular, the recollement (\ref{recollement diagram of D}) corresponds to the triangle TTF-triples $(\im(\bar{\psi}_{0}), \Ker(\bar{\psi}^{1}),\im(\bar{\psi}_{1}))$.
Observe that for any $0\leq i<p$ and $n\in\bbZ$, we have $\Ker(\bar{\psi}^{np+i})=\Ker(\bar{\psi^i})$ and  $\im(\bar{\psi}_{np+i})=\im(\bar{\psi_i})$. Hence there is an unbounded TTF-tuple which is periodic of period $p$:
$$(\cdots, \im(\bar{\psi}_{i-1}), \Ker(\bar{\psi}^{i}), \im(\bar{\psi}_{i}), \Ker(\bar{\psi}^{i+1}), \im(\bar{\psi}_{i+1}),\cdots ).$$
It follows that the ladder is periodic of period $p$.
\end{pf}

For any sequence $\underline{q}=(i_1,i_2,\cdots, i_q)$ of integers with
with $0\leq i_1<i_2<\cdots<i_{q}<p$, we can define
$\bar{\psi}^{\underline{q}}=\bar{\psi}^{i_1}\bar{\psi}^{i_2}\cdots \bar{\psi}^{i_q}$ and
$\bar{\psi}_{\underline{q}}=\bar{\psi}_{i_q}\cdots \bar{\psi}_{i_2}\bar{\psi}_{i_1}$ in a similar way as before.

The stable category of vector bundles over a weighted projective line of weight triple case is of particular interest, we refer to \cite{[KLM]} for more details. From now on, we assume $\bbY(p)$ is of weight triple, i.e, $\bbY$ has weight type $(p_1, p_2)$. In this case, we have the following result, which covers the main result in \cite{[C]}.

\begin{Thm}\label{main theorem for stab cat} Assume $\bbY$ has weight type $(p_1, p_2)$. For any $0\leq q<p$, let $\underline{q}=(1,2,\cdots, q)$ and $\underline{q}^c=(q+1,q+2,\cdots,p-1)$.
Then the following diagram is an infinite ladder of period $p$:
\begin{equation}\label{ladder diagram of stab cate}
\xymatrix@C=0.5cm{
  \underline{\vect}\bbY(q+1)
\ar@/_3pc/[rrr]_{\vdots}|{\bar{\psi}_{\underline{q}^c+\underline{1}}}
\ar@/^3pc/[rrr]^{\vdots}|{\bar{\psi}_{\underline{q}^c-\underline{1}}}
\ar[rrr]|{\bar{\psi}_{{\underline{q}}^c}}
&&& \underline{\vect}\bbY(p)
\ar@/_1.5pc/[lll]|{\bar{\psi}^{{\underline{q}}^c-{\underline{1}}} }
\ar@/^1.5pc/[lll]|{\bar{\psi}^{q^c}}
\ar[rrr]|{\bar{\psi}^{\underline{q}}}
\ar@/_3pc/[rrr]_{\vdots}|{\bar{\psi}^{{\underline{q}}+{\underline{1}}}}
\ar@/^3pc/[rrr]^{\vdots}|{\bar{\psi}^{{\underline{q}}-{\underline{1}}}}
&&& \underline{\vect}\bbY(p-q).
\ar@/_1.5pc/[lll]|{\bar{\psi}_{{\underline{q}}-{\underline{1}}}}
\ar@/^1.5pc/[lll]|{\bar{\psi}_{\underline{q}}}
}
\end{equation}
Consequently, for any $i\in\bbZ$, there is an exact sequence of triangulated categories
\begin{equation}\label{exact seq of stab cat}
\xymatrix{
\underline{\vect}\bbY(q+1)
\ar[rr]^{\bar{\psi}_{{\underline{q}}^c+i{\underline{1}}}} && \underline{\vect}\bbY(p) \ar[rr]^{\bar{\psi}^{{\underline{q}}+i{\underline{1}}}} & & \underline{\vect}\bbY(p-q)}.
\end{equation}
\end{Thm}

\begin{pf}
First we show that (\ref{exact seq of stab cat}) is exact for $i=0$. In fact, by Proposition \ref{adjoint tuples for stable category} we see that $\bar{\psi}_{{\underline{q}}^c}$ is fully faithful and $\bar{\psi}^{{\underline{q}}}$ is a quotient functor, since they are compositions of fully faithful functors and quotient functors respectively. Moreover, by similar arguments as in the proof of Lemma
\ref{kernel for bar of ri}, we can obtain that
$\Ker(\bar{\psi}^{{\underline{q}}})$ is the subcategory of $\underline{\vect}\bbY(p)$ consisting of objects of the form
$$\xymatrix{
E_0\ar[r]^{x_{0}} &E_1\ar[r]^{x_{1}} &  \cdots \ar[r]^{x_{q-1}} & E_{q}\ar[r]^{x_{q}} & E_{q+1}\ar@{=}[r] &  \cdots \ar@{=}[r] & E_{0}(x),
}$$
where $E_{0}$ is a direct sum of line bundles in $\coh\bbY$, which holds automatically by our assumption that $\bbY$ has weight type $(p_1, p_2)$, namely, each indecomposable vector bundle in $\coh\bbY$ is a line bundle. This yields an equivalence
$$\bar{\psi}_{q^c}: \quad\underline{\vect}\bbY(q+1)\longrightarrow \Ker(\bar{\psi}^{q}),$$ which proves the exactness of (\ref{exact seq of stab cat}) when $i=0$. Now by using Proposition \ref{adjoint tuples for stable category}, we obtain a recollement:
$$\xymatrix@C=0.5cm{
  \underline{\vect}\bbY(q+1)
\ar[rrr]|{\bar{\psi}_{{\underline{q}}^c}}
&&& \underline{\vect}\bbY(p)
\ar@/_1.5pc/[lll]|{\bar{\psi}^{{\underline{q}}^c-{\underline{1}}} }
\ar@/^1.5pc/[lll]|{\bar{\psi}^{q^c}}
\ar[rrr]|{\bar{\psi}^{\underline{q}}}
&&& \underline{\vect}\bbY(p-q).
\ar@/_1.5pc/[lll]|{\bar{\psi}_{{\underline{q}}-{\underline{1}}}}
\ar@/^1.5pc/[lll]|{\bar{\psi}_{\underline{q}}}
}$$
Since the category $\underline{\vect}\bbY(p)$ admits Serre duality, by using the reflecting approach in \cite{[J]}, we finally obtain a ladder, which should have the form (\ref{ladder diagram of stab cate}) by the adjointness sequence (\ref{adjoint sequence in stab}). Moreover, the periodicity of the ladder (\ref{ladder diagram of stab cate}) follows from the similar reason as before, namely, due to the facts
$\Ker(\bar{\psi}^{np+i})=\Ker(\bar{\psi^i})$ and $\im(\bar{\psi}_{np+i})=\im(\bar{\psi_i})$ for any $n\in\bbZ$.
We are done.
\end{pf}

This symmetric ladder above has direct applications in tilting theory. In fact, we can prove \cite[Theorem B]{[KLM]} in an inductive way. This also provides an effective way to construct new tilting objects in the stable category of vector bundles over a weighted projective line, due to the explicit expressions of the functors in the ladder (\ref{ladder diagram of stab cate}). These will be investigated in the forthcoming paper.

\section*{Acknowledgements}
The author is grateful to Henning Krause for leading me to this topic and also for his continuous discussions and valuable suggestions. The author also thanks Xiao-wu Chen, william Crawley-Boevey, Yan Han, Andrew Hubery and Helmut Lenzing for their useful discussions.

This work was partially supported by the Alexander von Humboldt Foundation in the framework of the Alexander von Humboldt Professorship endowed by the German Federal Ministry of Education and Research, and by the National Natural Science Foundation of China (No. 11801473) and Fundamental Research Funds for the Central Universities of China (No. 20720180006).

%%%%%%%%%%%%%%%%%%%%%%%%%%%%%%%%%%%%%%%%%%%%%%%%%%%%%%%%%%%%%%%%%%%%%%%%%%%%%%%%%%%%%%%%%%%%%%%%%%%%%%%
% References
%%%%%%%%%%%%%%%%%%%%%%%%%%%%%%%%%%%%%%%%%%%%%%%%%%%%%%%%%%%%%%%%%%%%%%%%%%%%%%%%%%%%%%%%%%%%%%%%%%%%%%%

\bibliographystyle{amsplain}

\end{document}